\documentclass[oupdraft]{bio}
\usepackage{url}
\usepackage{amsmath}
\usepackage{amsthm}

\newcommand{\argminD}{\arg\!\min}

\newcommand{\Var}{\operatorname{Var}}

\newcommand{\Cov}{\operatorname{Cov}}%

\usepackage{subcaption}
\usepackage{url}
\usepackage{tikz-cd}
\usetikzlibrary{positioning, arrows.meta}
\usepackage{booktabs}
\usepackage{url}

\newtheorem{definition}{Definition}[section]

\usepackage{algorithm}
\usepackage{algpseudocode}

\begin{document}
	
	\title{Unsupervised Domain Adaptation Under Hidden Confounding}
	
	\author{Carlos García Meixide$^{\ast,1,2}$, David Ríos Insua$^1$ \\[4pt]
		\textit{$^1$Instituto de Ciencias Matemáticas - CSIC, Madrid, Spain \\ 
			$^2$University of California, Berkeley, CA, US} \\[2pt]
		$^\ast$\texttt{meixide@berkeley.edu}
	}
	
	\maketitle

% Add a footnote for the corresponding author if one has been
% identified in the author list
%\footnotetext{To whom correspondence should be addressed.}

\begin{abstract}
{We introduce a new predictive mechanism that operates in the presence of hidden confounding across distributionally diverse data sources while ensuring consistent estimation of causal parameters—despite their recognized suboptimality for prediction in the literature. Our method is based on a novel estimand that captures the dependence structure between response noise and covariates, incorporating causal parameters into a generative model that adaptively replicates the conditional distribution of the test environment. Identifiability is achieved under a straightforward, empirically verifiable assumption. Our approach ensures probabilistic alignment with test distributions uniformly across arbitrary interventions, enabling valid predictions without requiring worst-case optimization or assumptions about the strength of perturbations at test time. Through extensive simulations, we demonstrate that our method outperforms state-of-the-art invariance-based and domain adaptation approaches. Additionally, we validate its practical applicability and superior target risk performance on a cardiovascular disease dataset.}{Causal Inference; Hidden Confounding; Domain Adaptation}
\end{abstract}

\section{Introduction}
\label{sec1}

Prediction in unseen domains presents challenges due to hidden variables that influence both the predictors and the response. Naïve empirical risk minimization on observed data produces estimates that can lead to completely incorrect predictions under distribution shift \citep{duchi,jeong2024outofdistribution,christiansen2021causal,saengkyongam2022exploiting,krueger2021out}, as it captures signals corresponding to the effects of hidden confounders. From a formal perspective, the challenge centers around the problem of \textit{prediction under hidden confounding}, which has surprsingly not yet been reconciled within the modern trend that interprets causal inference as a two-phase process. The first step, \textit{identification}, operates purely at the population level and does not involve data \citep{rubin74, robins86}. It aims at expressing estimands involving causal effects as a function of the probability distribution of the observed data. The second phase involves estimating this parameter using techniques such as targeted maximum likelihood estimation, Neyman-orthogonality, or debiasing \citep{laan06, syrg, chernodoubly, montanari, vansteelandt2023orthogonal}.  

In many practical scenarios, training data is collected from diverse sources, often exhibiting heterogeneous structures. Such data, fundamental to many data science problems, rarely originates from controlled experiments and instead contains unspecified perturbations \citep{Zhang2013}. While this may seem problematic, a key insight from the perturbation statistics literature is that i.i.d. data is inherently uninformative—perturbations are not obstacles to be mitigated but rather essential for achieving identifiability (see, for instance, \cite{ica}). Embracing this perspective, machine learning practitioners can effectively leverage data from varied environments, such as hospitals operating under different conditions, to enhance the robustness of predictive models \citep{didelez2006direct,he2008active}.

Traditional approaches to handling heterogeneous training distributions have followed two main strategies. The first assumes that the conditional distribution of the prediction target, given its causal parents, remains invariant across environments under arbitrary interventions on all variables except the response—an assumption known as invariance \citep{Haavelmo1943,Magliacane2017,Arjovsky2019,lin2022bayesian}. In \cite{peters2016causal}, a method was proposed to infer a full causal model, but it relied heavily on the absence of hidden confounding and encountered computational challenges as the dimensionality of covariates increased. In practice, methods based solely on invariance often perform similarly to empirical risk minimization (ERM), which ignores multi-environment information.

As a result, a second approach frames invariance as a regularized worst-case optimization problem, incorporating a hyperparameter that quantifies the expected degree of distribution shift in unseen domains relative to the training data \citep{henzinva,shen2023causality}. This concept, known as finite robustness, involves selecting a penalty term in a manner analogous to how $\ell_1$ regularization mitigates noise variance. A higher hyperparameter value shifts the solution closer to an instrumental variables (IV) estimator \citep{angrist1996identification}, which remains robust under arbitrarily large interventions but may be overly conservative when test predictors closely resemble the training data.

This article challenges the core claim of finite robustness, arguing that identifiability does not need to be sacrificed for accurate prediction. We propose a method that optimally adapts predictions to unlabeled test data while preserving identifiability under a general parameterization framework. This skeleton is identifiable under a mild and empirically verifiable condition (Theorem \ref{multidentif}). Our approach leverages invariance across multi-source training data for both identification and estimation, while unlabeled test data informs a generative model to produce predictions. These two components define our method, Generative Invariance (GI). Although invariance remains central to our approach, we show that the key invariant property is not the conditional distribution given the source and covariates, but rather smaller, fundamental components of this distribution. 

Our approach fundamentally departs from the conventional view of predictive robustness as a minimax problem, demonstrating that there is no need to introduce bias into causal parameter estimators for accurate prediction. We achieve this through a novel component that mitigates the inherent overconservatism of classical instrumental variable (IV) methods while avoiding the need for finite robustness-style regularization. In this sense, Generative Invariance (GI) serves as an unsupervised domain adaptation method, leveraging probabilistic information from unlabeled observations in the target domain. This is particularly important in practical applications, where predictions are typically required for multiple units rather than a single fixed instance.

In summary, our main contributions are as follows: 

\begin{itemize}
	\item We introduce a generative model for endogenous residuals that enables optimal predictions in terms of test mean squared error (MSE), relying only on classical instrumental variable (IV) conditions. This represents a fundamental shift in perspective, as causal parameters have traditionally been viewed as suboptimal for prediction under test distribution shifts \citep{rothenhausler2021anchor}.
	\item Unlike classical distributional robustness approaches, GI does not require partial knowledge of the direction or magnitude of test distribution shifts.
	\item Our method does not assume additive shifts. Unlike recent approaches \citep{kostin}, our framework does not impose an internal measurement error structure within the covariates.
	\item In multi-environment settings with heterogeneous training sources, GI treats all data sources symmetrically. In contrast to methods that require a reference environment with a dominant second moment (e.g., \cite{shen2023causality}), our approach imposes no such restriction.
	\item Our distributional setting extends beyond covariate shift: due to unobserved confounding, the conditional distribution also varies across environments.
	\item GI achieves identifiability in the one-dimensional case without requiring multi-source data and does not rely on hyperparameter tuning in multivariate setups.
	\item Empirically, we demonstrate that GI consistently outperforms state-of-the-art methods across hyperparameter values \citep{shen2023causality} and achieves superior results compared to leading unsupervised domain adaptation techniques.
	\end{itemize}

\subsection{Structure of the paper}

Section 2 establishes the core assumptions underlying our methodology. Section 3 introduces a new estimand for optimal prediction under distribution shifts, covering its theoretical foundations and detailing how our generative model enables the replication of the test environment's conditional distribution.  Section 4 presents plug-in estimators leading to empirical GI, implementation details and closed-form solutions for the proposed estimator. Large-sample properties are covered in Section 5, providing a key finite-sample concentration inequality leading to an asymptotic normality result.  Section 6 presents extensive simulations showcasing the superior predictive accuracy of our approach with respect to state-of-the-art methods as well as a medical application. Minor proofs are provided in the Supplement, while those concerning asymptotic theory are included in Appendices A and B of this document.

\section{Notation and setup}\label{sec:not}

In exploring the underlying mechanisms of random variables $X \in \mathbb{R}^p$ and $Y \in \mathbb{R}$, consider the following general structural causal model (SCM): 
\begin{equation} \label{firstscm}
	X \leftarrow  (1-\theta)X_0 +  \theta X_1 \textrm{ and } Y \leftarrow  \beta^T_*X + \varepsilon_Y,
\end{equation}

\noindent where $\theta \in \{0,1\}$ and $X_0, X_1 \in \mathbb{R}^p$ are mutually independent random variables, with $1$ and $0$ respectively representing training and testing phases. The response residual $\varepsilon_Y \in \mathbb{R}$ is independent of $\theta$  but might not be independent of $X_0$ or $X_1$.    Conditional on $\theta$, we permit dependencies among the components of $(X^1, \ldots, X^{p}, \varepsilon_Y)$ unlike the standard approach (i.e. absence of hidden parent variables). Noise endogeneity between $X$ and $Y$ is equivalent to the existence of a hidden confounder that makes residual contributions not assumable to be independent, allowing us to omit hidden variables from a formal viewpoint. The error term $\varepsilon_Y$ is zero-mean. Importantly, note that $X_0$ and $X_1$ are not to be interpreted as \textit{potential outcomes} \citep{splawa} as they do not try to encode counterfactual duality for a single individual: they should be rather thought of as independent random variables so that sample units are not related across domains.

Our objective is to generate predictions for given covariate values \(X_{0,1}, \ldots, X_{0,n}\) in the test phase (\(\theta = 0\)), using observations \((X_{1,1},Y_{1,1}), \ldots, (X_{1,n},Y_{1,n})\) collected during the training phase (\(\theta = 1\)). Notably, our framework is designed to satisfy the classical instrumental variable (IV) assumptions \citep{wright1928tariff,angrist1996identification}, with \(\theta\) serving as an instrument. Specifically, \(\theta\) (i) is associated with \(X\), (ii) is independent of unmeasured confounders, and (iii) has no direct effect on \(Y\).  

The IV approach is fundamental in addressing hidden confounding in observational studies. It relies on the instrument \(\theta\) being orthogonal to the residuals of the response variable \(Y\) obtained after adjusting for the covariates \(X\) through the \textit{causal parameter} \(\beta_*\). This parameter is termed \textit{causal} because it directly determines how the covariates influence \(Y\) in the structural causal model (SCM) (\ref{firstscm}). However, a standard multivariate linear regression model would yield an inconsistent estimate of \(\beta_*\) due to endogeneity arising from target noise.  

According to IV theory, treating the test data as a distinct environment would be sufficient to identify the causal parameter \(\beta_*\). However, a key challenge is that we do not observe the response values \(Y_0\) in a sample drawn from the test distribution. Despite this limitation, we will demonstrate that \(\beta_*\) remains identifiable when \(p=1\), even in the absence of test labels.  

%Under the three classical IV assumptions above, in the simpler one-dimensional setting the goal is to minimize with respect to $\beta \in \mathbb{R}$ the risk $\mathbb{E}\left[\mathbb{E}^2\left[Y-\beta X \mid \theta\right]\right].$ This looks for the $\beta$ that hinders $\theta$ the most when trying to explain the variability of residuals. This is because conditional expectations can be understood as a projection in an  
%$L^2$ space of random variables, representing the best 
%approximation of a random variable given the conditioning one in terms of expected squared distance. The technical nuance that the square is taken outside the inner expectation elucidates the rationale behind the feasibility of consistently estimating the causal parameter. From a formal viewpoint, this subtlety is pivotal to \textit{identify} $\beta_*$: the interested reader can easily verify this by assuming, for simplicity, that $\mathbb{P}(\theta=0)=\mathbb{P}(\theta=1)=\frac{1}{2}$ and then expanding the expression inside the expectation and squaring it. This exercise helps to bear in mind the often forgot important fact that the manifestation of identifiability happens at population level. 

% where $\varepsilon_X:=  H\alpha_* + \varepsilon_X $ and $\varepsilon_Y:= H\gamma_* + \varepsilon_Y$ with a slight notational liberty. It is important to note that the newly defined residuals are not uncorrelated anymore: $$\Cov(\varepsilon_X,\varepsilon_Y) = \Cov(H\alpha_* + \varepsilon'_X,H\gamma_* + \varepsilon'_Y) = \gamma_*\alpha_* \Var H =: K_*,$$ accounting for hidden confounding. 

{Denote by $\mathbb{P}$ the joint distribution of $\theta$, $X$ and $Y$ and} let $(X_z,Y_z)$ be the random variable following the conditional distribution of $(X,Y)$ given $\theta=z $. We will notationally switch from $\mathbb{E}[g(X,Y) | \theta=z]=\mathbb{E}g(X_z,Y_z)$ to $\mathbb{E}_zg(X,Y)$ deliberately for any function $g$ that is measurable with respect to the appropriate $\sigma$-algebras whenever the previous expectations exist. Denote the joint distribution of $(X_0,Y_0)$ and $(X_1,Y_1)$, respectively,  by $\mathbb{P}^0$ and $\mathbb{P}^1$. This allows to decouple the initial SCM (\ref{firstscm}) into a pair of SCMs representing the test and train distributions respectively: $\mathbb{P}^0: Y_0 \leftarrow \beta^T_*X_0 + \varepsilon_{Y,0}, \quad \mathbb{P}^1: Y_1 \leftarrow \beta^T_*X_1 + \varepsilon_{Y,1}.$ Note that we tacitly assume that there are no interventions on the residuals distribution at test time, although this can be relaxed as explained in Section \ref{multiv}. We  assume that confounding is \textit{invariant} across the training and test phases, that is, $$\Cov (X_0,\varepsilon_{Y,0}) =\Cov (X_1,\varepsilon_{Y,1}) =: K^*\in \mathbb{R}^p, \quad \mathbf{(A1)}$$

\noindent where the $j$-th entry of $\Cov (X_z,\varepsilon_Y)$ is $\mathbb{E}X^j_z\varepsilon_Y$. Assumption $\mathbf{(A1)}$ is referred to as \textit{inner-product invariance} in \cite{rothenhausler2019causal}.  Observe that $\Cov(X,\varepsilon_Y \mid \theta=z) =(1-z) \Cov(X_0,\varepsilon_{Y,0}) + z \Cov(X_1,\varepsilon_{Y,1})  = (1-z)K_* + zK_* =K_*$ for $z=0,1$. Note that the additive interventions model \citep{rothenhausler2019causal} is a particular case of our SCM (\ref{firstscm}) plus assumption (A1). Namely, let  $X_1\leftarrow \varepsilon_X$ and $X_0 \leftarrow \varepsilon_X + V$ with $V$ independent of $\varepsilon_X$ and  $\varepsilon_Y$; then $\Cov(\varepsilon_{Y,0},X_0)= \Cov(\varepsilon_{Y,0},\varepsilon_X + V) = \Cov(\varepsilon_{Y,0},\varepsilon_X )  = \Cov(\varepsilon_{Y,1},X_1)$. 

While estimators whose probability limit is the causal parameter \(\beta_*\) are useful for understanding underlying mechanisms, relying solely on them for predictive purposes is suboptimal—particularly when test-phase perturbations are mild \citep{rothenhausler2021anchor}. As a example, suppose that the train and test distributions coincide ($\mathbb{P}^1=\mathbb{P}^0=:\mathbb{P}$); although the ordinary least squares (OLS) estimator is biased because of the assumption $\varepsilon_Y \perp X$ being violated, by definition we have $\mathbb{E}\left[\left(Y- \beta^T_*X\right)^2\right] \geq \mathbb{E}_0\left[\left(Y- \beta^T_{\operatorname{OLS}}X\right)^2\right], $ where $\beta_{\operatorname{OLS}}:=\argminD_{\beta} \mathbb{E}_1\left[(Y-\beta^TX )^2\right]=(\mathbb{E}X_1X_1^T)^{-1}\mathbb{E}X_1Y_1 \stackrel{K_*=0}{=} \beta_*.$ This exemplifies a scenario in which relying on the causal parameter results in suboptimal predictive performance, see Figure \ref{fig:simpson}. However, OLS yields an arbitrarily high MSE when the test distribution diverges substantially from the training one. Established methods such as \textit{anchor regression} \citep{rothenhausler2021anchor} intentionally introduce bias to mitigate the excessive conservativeness that arises from relying solely on the causal parameter for prediction, particularly in test distributions that closely resemble the training distribution. This exploits the fact that endogeneity bias is implicitly and beneficially leveraged by OLS in test environments where the probabilistic structure remains similar to that of the training phase.

We argue that the strategy of deliberately distributing prediction error across perturbations to avoid poor performance in specific scenarios is no longer necessary.  An illustration of the behaviour of our method is in Figure \ref{firstfig}. Notice therein how the more different the test covariates are with respect to the training set, the worse OLS performs. Our procedure (green) reproduces the probabilistic structure of the test environments by looking at the blue data and just at the covariate values of each test sample:  {GI moves away from traditional reliance on regularization and \textit{finite robustness} to implement optimal predictions in an unlabelled domain without sacrificing consistent estimation of the causal parameters, as shown next. }

\section{Population Generative Invariance}  \label{popu}
This section frames optimal prediction under distribution shifts through a generative model which we parametrize with estimands recoverable as minimizers of a certain theoretical risk. We study conditions under which these estimands are identifiable. 

Let $\Sigma_1:=\Cov X_1$ be the $p\times p$ covariance matrix of the covariates at training. In general, $\Sigma_1$ will differ from $\Sigma_0:=\Cov X_0$. Assume throughout this section that there are no interventions on the distribution of the response noise variance across test and train, a hypothesis to be relaxed in Remark \ref{diff_var}, and let $\Var \varepsilon_Y:=  \sigma^2_Y$. We shall construct a new matrix, the joint covariance matrix of the $p+1$ dimensional vector $(X_1,Y_1)$, augmenting $\Sigma_1$ with a real $p$-dimensional vector $K$ such that the matrix $
\left(\begin{array}{cc}
	\Sigma_1 & K \\
	K^T & \sigma_Y^2
\end{array}\right)
$ remains positive definite (p.d.). To ensure it, consider the Schur complement of $\sigma^2_Y$ in the new matrix, $\mathcal{S}=  \Sigma_1-\frac{K K^T}{\sigma_Y^2}
$. As $\Sigma_1$ is p.d., if $\mathcal{S}$ is p.d. then the augmented matrix will be p.d. as well \citep{john}. Thus, we just need to ensure that $\mathcal{S}$ is also p.d., which means that for all non-zero vectors $ {v} \in \mathbb{R}^p$, $
{v}^T   \Sigma_1  {v}-\frac{1}{\sigma_Y^2}\left( {v}^T K\right)^2>0
$ holds. This is equivalent to what we call
\begin{equation}\label{schur_cond}
	\textrm{``causal Mahalanobis condition 1'':   }K^T   \Sigma_1^{-1} K<\sigma_Y^2
\end{equation}

\noindent Thus, $K$ must lie in the ellipsoid of $\mathbb{R}^p$ defined by restriction (\ref{schur_cond}), denoted by $\mathcal{K}_1$ onwards so that $\mathcal{K}_z := \{K \in \mathbb{R}^p: K \textrm{ satisfies causal Mahalanobis condition }z\}$ for $ z \in \{0,1\}$.  Let $\mathcal{K}_{01}:=\mathcal{K}_0\cap \mathcal{K}_1$- which is non-empty because of (A1)- and define, for each $K \in \mathcal{K}_{01}$, the \textit{generator}
\begin{equation}\label{our_gen}
	g_K: (x,\xi) \in \mathbb{R}^{p} \times \mathbb{R} \longmapsto K^T  \Sigma_0^{-1} (x - \mathbb{E}X_ 0)+ \xi \sqrt{\sigma^2_Y - K^T \Sigma^{-1}_0K} \in \mathbb{R}
\end{equation}

% $$g_K: (x,\xi) \in \mathbb{R}^{p+1} \longmapsto (x - \mathbb{E}X_1)^T   \Sigma_1^{-1} K+ \xi S \in \mathbb{R}$$

%\noindent where $S^p_K$ is a positive number that depends on the covariates dimension $p$ and the value of $K$ indexing $g_K$.  $S^p_K$ ensures the variance of $g_K\left(X_ 0, \xi\right)$ to be the same as that of $\varepsilon_Y$. It turns out that $$S^1_K= \sqrt{\frac{\sigma^2_Y \Var X_ 0 - K^2}{\Var X_ 0}}$$ 
%$$S^p_K= \sqrt{\sigma_Y - K^T \Sigma^{-1}_0K}$$
%$$S^2_K= \sqrt{\frac{V_1V_2\sigma^2_Y - \sigma^2_YL^2 - V_2K_1^2 - V_1K_2^2 + 2K_1K_2L}{V_1V_2-L^2}}, \quad \textrm{with} \left(\begin{array}{cc}
	%	V_1 & L \\
	%	L & V_2
	%\end{array}\right)=\Sigma_ 0$$
	
	%\noindent Note that for $K \in \mathcal{K}_0$, the radicand above is always positive. Although they become intricate for $p \geq 3$, the explicit expressions of $S^p_K$ are only of theoretical interest. This is because users are just interested on the conditional expectation of the generator $g_K$ given $X_ 0=x$ for practical applications. For the interested reader, see further details in Appendix \ref{sigmay}. 
	
	%is a constant whose form depends on $p$:
	
	% \begin{itemize}
		% \item For $p=1$, $$S_p= \sqrt{\frac{V_Y V_1 - K_1^2}{V_1}}$$
		% \item For $p=2$, $$S_p= \sqrt{\frac{V_1V_2V_Y - V_YL^2 - V_2K_1^2 - V_1K_2^2 + 2K_1K_2L}{V_1V_2-L^2}}$$
		% \end{itemize}
	
	% \begin{remark}
		%
		% \end{remark}
	\noindent Proposition \ref{gk} asserts that evaluations of these functions on $X_ 0$ with zero mean, unit variance noise $\xi$ are intimately connected to the joint probabilistic structure of the covariate and response residuals at test phase ($z=0$).
	
	\begin{Proposition} \label{gk}
		Let $\xi$ be a real random variable such that $\mathbb{E}\xi =0$ and $\Var \xi = 1$. For $K \in \mathcal{K}_{01}$, $g_K(X_ 0,\xi)$ is a random variable such that i) $\mathbb{E}g_K(X_ 0,\xi)=0$; ii) $\operatorname{Var}g_K(X_ 0,\xi) = \sigma^2_Y$; iii) $\Cov(X_ 0,g_K(X_ 0,\xi )) = K$.
	\end{Proposition} 
	
	\noindent Proving i) and ii) is inmediate; iii) holds because $\mathbb{E}(X_ 0 - \mathbb{E}X_ 0)K^T  \Sigma_ 0^{-1} (X_ 0 - \mathbb{E}X_ 0) = \mathbb{E}(X_ 0 - \mathbb{E}X_ 0) (X_ 0 - \mathbb{E}X_ 0)^T  \Sigma_ 0^{-1}K = \Sigma_ 0\Sigma_ 0^{-1}K= K$. Proposition \ref{gk} succintly ensures that evaluations of the generator $g_{K_*}$ resemble the original endogenous residuals in terms of moments up to second order. As a consequence, assuming Gaussianity and setting $\beta=\beta_*, K=K_*$, we recover a copy of $(X_ 0,Y_ 0)$ as we see next. 
	\begin{Corollary} \label{copy} Let $\xi \sim \mathcal{N}(0,1)$. If $\mathbb{P}^0$ is Gaussian then $(X_ 0,g_{K_*}(X_ 0,\xi) + \beta^T_*X_ 0) \sim \mathbb{P}^ 0. $
	\end{Corollary}
	\noindent  {This result motivates the name \textit{Generative Invariance}. Explicitly quoting \cite{goodfellow2014generative}, ``\textit{we explore the special case when the generative model generates samples by passing random noise through a multilayer perceptron}''. In Corollary \ref{copy}, the ``\textit{prior on input noise variables}'' is $\mathcal{N}(0,1)$, while the ``\textit{mapping to data space}'' is $g_K(X,\cdot)$, with $K \in \mathcal{K}_{01}$ . Our generators follow a much simpler linear dependence structure on $K$ and, importantly, they are \textit{conditional} generative models because they also take as input the covariates distributed at test. Our generative model also respects the \textit{causal generative process from exogenous noise} involved in the \textit{Principle of Independent Components} \citep{peters2017elements}}. 
	Let $
	\mathcal{M}= \{m_{\gamma}: (x,\xi) \in \mathbb{R}^{p+1} \longmapsto g_K(x,\xi) + \beta^Tx \textrm{ such that } \gamma=(\beta, K) \in \Xi := \mathbb{R}^p\times \mathcal{K}_{01}\}
	$. Each element in $\mathcal{M}$ is a function $m_\gamma$  indexed by the tuple $\gamma=(\beta,K) \in \Xi$ taking covariates $x$ and noise $\xi$, inducing a conditional distribution as follows. Fix $\gamma$, draw $\xi \sim \mathcal{N}(0,1)$ and consider the conditional distribution of $m_{\gamma}(X_ 0,\xi)$ given $X_ 0=x$. Each $m_{\gamma}$ represents a procedure involving a generative model in two phases: first sample from the distribution of $X_0$ and then sample $Y_0$ from the conditional distribution given $X_0$ induced by the generator. If observations $X_{0,1}, \ldots, X_{0,n}$ from the marginal covariate distribution of the test stage and independent $\xi_1, \ldots, \xi_n \sim \mathcal{N}(0,1)$ are generated, then $(X_{0,1},m_{\gamma}(X_{0,1},\xi_1)),\ldots,(X_{0,n},m_{\gamma}(X_{0,n},\xi_n))$ is distributed according to the test distribution $\mathbb{P}^0$ when $\gamma=(\beta_*,K_*)$. 
	
	% ; see for example \citep{goodgan,shen2023engression}. 

	\begin{remark}\label{doint}
		The do-interventional distribution \citep{Pearl_2009} of $Y_0$ under $\operatorname{do}(X_0=x)$ is the same as the distribution of $m_{(\beta_*,K)}(x,\xi)$  {setting $K \equiv 0$}, $\xi \sim \mathcal{N}(0,1)$.  {Considering a non-zero $K \in \mathcal{K}_{01}$ is precisely what will distinguish our approach from standard IV prediction.}
	\end{remark}

	\subsection{One dimensional feature ($p=1$)}
	
	In addressing the estimation of model parameters, we show that their actual value can be identified through a theoretical risk minimization problem, with attainable and unique solution under lean assumptions as Theorem \ref{t1} will show. Provided it exists, we define population GI using the parameters
	\begin{equation}\label{thrisk}
		(\beta_{opt}, K_{opt}):=\argminD_{(\beta,K) \in \Xi} \mathbb{E}_1\left[(Y - \beta X - K(X - \mathbb{E}X ))^2\right].
	\end{equation}
	
	\noindent Expanding the square inside the expectation and setting derivatives to zero we arrive at the following result.
	
	\begin{Proposition}\label{p1} Whenever it exists, $(\beta_{opt}, K_{opt})$ must satisfy $
		(\beta_* - \beta_{opt})\mathbb{E} ^2X_1 = 0
		$ and $K_{opt}{\Var X_1} -  K_*= {(\beta_* - \beta_{opt})\Var X_1 }
		$

	\end{Proposition}
	\noindent As a consequence, identifiability of $(\beta_*,K_*)$ is inmediate provided that covariates at traning are not zero-mean.
	\begin{theorem}\label{t1}
		If $ \phantom{s} \mathbb{E}X_1\neq 0$ then $(\beta_{opt}, K_{opt})$ exists, is unique, $\beta_{opt}= \beta_* \textrm{ and }K_{opt} \Var X_1= K_*.$
		
	\end{theorem}
	\noindent Sampling from the covariate distribution in the test environment is what facilitates the creation of artificially correlated residuals and, therefore, optimal predictions. This is the reason why many methodologies need  ``causal regularization'' when making predictions for a single $x \in \mathbb{R}$ paying the price of identifiability bias, see e.g. \cite{henzinva}.
	
	\subsection{Multivariate features ($p>1$)}\label{multiv}

	Suppose that at training there is access to several data sources such that we observe samples from both covariates and responses in all of them. Let us show that the identification of $\beta_*$ and $K_*$ together with their subsequent estimation are still possible for {covariate dimension} $p>1$ under purely algebraic hypotheses. For this, let us extend our original setup $\theta \in \{0,1\}$ to $\theta \in \{0,1,2,\ldots, Z\}$, where $Z$ is the number of different training environmnents and $\theta=0$ is the test environment. The original SCM (\ref{firstscm}) becomes 
	$$X \leftarrow \sum_z 1\{\theta = z\} X_z\text { and } Y \leftarrow \beta_*^T X+\varepsilon_Y, \text{ with } \mathbb{E}\varepsilon_Y = 0  \text{ and } \mathbb{E}\left[ X \varepsilon_Y \mid \theta = z \right] = K_*$$
	
	\noindent 	Now the \textit{causal Mahalanobis condition} (\ref{schur_cond}) must be satisfied for every $z=0,1,\ldots,Z$ and, therefore, $K_*$ lives in the intersection of all the ellipsoids involved, that is $K_* \in \bigcap_{z=0}^{Z} \mathcal{K}_z $. In the subsequent writing, we implicitly condition on the event $ \bigcup_{z=1,\ldots,Z}\{\theta=z\} $ to exclude the test phase $\theta=0$.  Assume for simplicity that each training observation is generated equally likely by each environment $1, \ldots, Z$ and thus $\mathbb{P}(\theta=z)=\frac{1}{Z}$ for $z=1,\ldots,Z$. Section \ref{remnz} clarifies why we can assume this without loss of generality and how to adapt the estimators to contexts with different number of observations per environment. The theoretical risk minimization problem (\ref{thrisk}) becomes
	\begin{equation}\label{multhrisk}
		(\beta_{opt}, K_{opt}) = 	\arg \min _{(\beta, K) \in \Xi} \sum_{z=1}^Z \mathbb{E}\left[\left(Y_z-\beta^T X_z-K^T\left(X_z-\mathbb{E} X_z\right)\right)^2\right].
	\end{equation}
	
	\noindent The presence of the term $\mathbb{E} X_z$ in (\ref{multhrisk}) key. If $\mathbb{E}X$ was there instead, the objective function would coincide with an expression that has the same form as $Z$ times the univariate risk (\ref{thrisk}). To see this, use the iterated expectations lemma $$\mathbb{E}\mathbb{E}\left[\left(Y - \beta^TX - K^T(X - \mathbb{E}X )\right)^2 \mid \theta\right] =  \frac{1}{Z}\sum_{z=1}^Z\mathbb{E}\left[\left(Y_z - \beta^TX_z - K^T(X_z - \mathbb{E}X ) \right)^2\mid \theta=z\right] $$
	
	\noindent  {Our goal is to estimate the invariant, causal parameters $\beta_{*}, K_* \in \mathbb{R}^p$ using labelled data (with observed $Y$) in the training sources $\{1, \ldots, Z\}$ and launch optimal predictions making use of these estimates and of unlabelled data in the testing source $\{0\}$ through the generative model \ref{our_gen}}. 	The next result is important to study conditions sufficient for identifiability of $(\beta_*,K_*)$ in the multivariate case, analogously to Proposition \ref{p1}. 
	\begin{Proposition}\label{pop_mult}
		Whenever it exists, the minimizer of problem (\ref{multhrisk}) must satisfy \\ $\sum_{z=1}^Z \mu_z \mu_z^T\left(\beta_*-\beta_{opt}\right)=0$ and $\phantom{s} \sum_{z=1}^Z \Sigma_z K_{opt}-ZK_*=\sum_{z=1}^Z \Sigma_z\left(\beta_*-\beta_{opt}\right).$
	\end{Proposition}
	
	\noindent We recover Proposition \ref{p1} upon setting $Z=p=1$ in the previous result. Next, by direct inspection of the first equation in Proposition \ref{pop_mult}, asking $\sum_{z=1}^Z \mu_z \mu_z^T$ to have an empty null space is sufficient to force $\beta_*=\beta_{opt}$ and, therefore, ensure identifiability of $K_*$ as well.
	\begin{Corollary}\label{howkstar}
		If the matrix $\sum_{z=1}^Z \mu_z \mu_z^T$ is full rank, there exists a unique solution $(\beta_{opt},K_{opt})$ to problem (\ref{multhrisk}) such that $\beta_{opt}=\beta_*$ and $\sum_{z=1}^Z \Sigma_z K_{o p t}=Z K_*.$
	\end{Corollary}

	\noindent For $Z>1$, recall that $\frac{1}{Z}\sum_{z=1}^Z \Sigma_z$ is in general not equal to the empirical covariance matrix of the covariate observations pooled across environments.  
	
	\begin{remark}\label{diff_var}
		
		If interventions on the response noise variance were allowed, during the proof of Proposition \ref{pop_mult}, the term $Z \sigma^2$ would be replaced by $\sum_{z=1}^Z \sigma^2_z$, where $\sigma^2_z:= \Var \varepsilon_{Y,z}$ is the response noise variance in source $z$. This means that identifiability of $\beta_*$ and $K_*$ would be kept free lunch despite sacrificing the assumption that response noise variance is the same across environments. 
	\end{remark}
	
	A sufficient condition for $\sum_{z=1}^Z \mu_z \mu_z^T$ to be invertible is that $\operatorname{dim}\langle \mu_1, \ldots, \mu_Z \rangle = p$, ensured by new linear algebra results in the Supplement. Therefore, as an inmediate consequence of Proposition \ref{pop_mult}, we have:

	\begin{theorem} \label{multidentif}
		If the population mean vectors  $\mu_1,  \ldots, \mu_Z$ of the covariates corresponding to the different environments span $\mathbb{R}^p$, $\beta_*$ and $K_*$ are identifiable. 
	\end{theorem}

	\noindent This is reminiscent of the well-known result in the IV literature asserting that if the dimension of the instrument is greater than or equal to the dimension of the covariates, then the causal parameter is the probability limit of two-stage estimators \citep{hall2005generalized}.  {In addition, Theorem \ref{multidentif} is a {constructive} identifiability result: we explicitly provide a sufficient condition that can be checked given just observables.} {If the hypothesis in Theorem \ref{multidentif} does not hold, potentially due to having access to $Z<p$ environments, then we are in a \textit{weak instrument} scenario.
		
		To conclude this section, we set $p=1$ in Theorem \ref{multidentif} and reach the conclusion that GI identifiability is possible with a single data source.

		\begin{Corollary}
			In a univariate model without intercept, $\beta_*$ and $K_*$ are identifiable using data from just one environment such that the population mean of its covariate is nonzero. 
		\end{Corollary}

		\section{Empirical Generative Invariance}
		
		Moving forward to inference aspects, we rewrite problem (\ref{multhrisk}) as:
		$$\begin{aligned}&\arg \min _{(\beta, K) \in \Xi}\sum_{z=1}^Z \mathbb{E}\left[\left(Y_z-\beta^T X_z-K^T\left(X_z-\mathbb{E} X_z\right)\right)^2\right] =\\
			&\arg \min _{(\beta, K) \in \Xi}\sum_{z=1}^Z\int \left(Y-\beta^T X-K^T\left(X- \int Xd\mathbb{P}^z\right)\right)^2 d\mathbb{P}^z \\&=: \arg \min _{(\beta, K) \in \Xi}T(\mathbb{P}^1, \ldots, \mathbb{P}^Z) .\end{aligned}$$
		
		\noindent Assuming that $T$ is well-defined at the conditional empirical measures $\mathbb{P}_{n}^z \textrm{ for }z=1,\ldots, Z$, consider the \textit{plug-in} \citep{vandegeer2010empirical} estimator $\arg \min _{\beta, K} T(\mathbb{P}_{n}^1, \ldots, \mathbb{P}_{n}^Z)$, which leverages samples \begin{align}\label{data}
			(X_{1,1},Y_{1,1}),\ldots,(X_{1,n},Y_{1,n}) &\sim \mathbb{P}^1 \notag \\
			& \vdots \notag \\
			(X_{Z,1},Y_{Z,1}),\ldots,(X_{Z,n},Y_{Z,n}) &\sim \mathbb{P}^Z
		\end{align}
		
		\noindent assuming we have $n$ observations per environment without loss of generality, see Remark \ref{remnz}. Let us concatenate the whole data in $X \in \mathbb{R}^{nZ \times p}$ and $Y \in \mathbb{R}^{nZ}$. Define \begin{equation}\label{eich}
			H_Z=\frac{1}{n}I_Z \otimes 1_n1^T_n  \textrm{ and } M_Z=H_Z X,\end{equation} 
		%	
		%	\noindent For example, if $Z=2$ 
		%	
		%	\[H_2=
		%	\frac{1}{n}
		%	\begin{pmatrix}
			%		1 &  \\
			%		& 1 
			%		
			%	\end{pmatrix}
		%	\otimes
		%	\begin{pmatrix}
			%		1 & 1 & 1 & \cdots & 1 \\
			%		1 & 1 & 1 & \cdots & 1 \\
			%		1 & 1 & 1 & \cdots & 1 \\
			%		\vdots & \vdots & \vdots & \ddots & \vdots \\
			%		1 & 1 & 1 & \cdots & 1
			%	\end{pmatrix} = \frac{1}{n}
		%	\begin{pmatrix}
			%		1 & 1 & 1 & \cdots & 1 & & & & & \\
			%		1 & 1 & 1 & \cdots & 1 & & & & & \\
			%		1 & 1 & 1 & \cdots & 1 & & & & & \\
			%		\vdots & \vdots & \vdots & \ddots & \vdots & & & & & \\
			%		1 & 1 & 1 & \cdots & 1 & & & & & \\
			%		& & & & & 1 & 1 & 1 & \cdots & 1 \\
			%		& & & & & 1 & 1 & 1 & \cdots & 1 \\
			%		& & & & & 1 & 1 & 1 & \cdots & 1 \\
			%		& & & & & \vdots & \vdots & \vdots & \ddots & \vdots \\
			%		& & & & & 1 & 1 & 1 & \cdots & 1 \\
			%	\end{pmatrix}
		%	\]
		%	
		\noindent where $\otimes$ designates the Kronecker product of two matrices. Using these definitions, 
		\begin{equation}\label{emprisk}
			(\widehat{\beta}, \widehat{K_{opt}}):= \arg \min _{\beta, K} nT\left(P_n^1, \ldots, P_n^Z\right)=\argminD_{\beta,K} \lVert Y - X\beta - (X - M_Z )K\rVert_2 ^2\end{equation}
		
		\noindent and we take $\widehat{\beta}$ and $\widehat{K}:= \frac{1}{Z} \sum_{z=1}^Z\widehat{\Sigma}_z \widehat{K_{opt}}$ as our estimators for $\beta_*$  and $K_*$ following the population expression in Corollary \ref{howkstar}. Let $\widehat{\mu}_z\in \mathbb{R}^p$ be the empirical average of the covariates in environment $z=1,\ldots,Z$. We have
		
		\begin{theorem}[Properties of $H_Z$ and $M_Z$]\label{hprop} Let $H_Z$ and $M_Z$ be as in (\ref{eich}). Then: i) $H_Z$ is idempotent; ii) $X^TM_Z=M_Z^T M_Z$; iii) $X^TX - M^T_ZM_Z= n\sum_{z=1}^Z\widehat{\Sigma}_z$; iv) $n \sum_{z=1}^Z \widehat{\mu}_z\widehat{\mu}^T_z= M^T_ZM_Z$
		\end{theorem}
		\begin{proof} \textit{i.} $H_Z^2=H_Z^TH_Z=(\frac{1}{n}I_Z \otimes 1_n1^T_n)(\frac{1}{n}I_Z \otimes 1_n1^T_n) = \frac{1}{n^2}I_Z \otimes n 1_n1^T_n = \frac{1}{n}I_Z \otimes 1_n1^T_n = H_Z$. The proof for the remaining items involves straightforward algebraic manipulations.
		\end{proof}

		\begin{remark}\label{remnz} In case there is a different number of observations per source $n_z$, it is easy to verify that iii) in Theorem \ref{hprop} becomes $X^TX - M^T_ZM_Z = \sum_{z=1}^Z n_z \widehat{\Sigma}_z$
			where $\widehat{\Sigma}_z$ is computed with denominator $n_z$ (rather than $n_z-1$). This is because, at the population level and under identifiability, $K_*=\left(\sum_{z=1}^Z w_z \Sigma_z \right)K_{o p t}$ with $w_z=\mathbb{P}(Z=z)$. Empirically, we estimate these weights as $ \widehat{w}_z=\frac{n_z}{ \sum_{z=1}^Z n_z}$ so that $\widehat{K}=\frac{1}{\sum_{z=1}^Z n_z}(X^TX - M_Z^TM_Z)\widehat{K}_{opt}.$ Analogously, $\sum_{z=1}^Z n_z
			\widehat{\mu}_z \widehat{\mu}_z^T=M_Z^T M_Z$. Regarding problem (\ref{thrisk}), its natural generalization to heterogeneous sample sizes would be a weighted version of \ref{multhrisk}. When replacing each conditional expectation with environment-wise average and plugging  $\widehat{w}_z$, the factor $\frac{1}{n_z}$ from the empirical averages will cancel out and therefore (\ref{emprisk}) is correct for different sample sizes. 
		\end{remark}
		%$\arg \min _{(\beta, K) \in \Xi} \sum_{z=1}^Z w_z \mathbb{E}\left[\left(Y_z-\beta^T X_z-K^T\left(X_z-\mathbb{E} X_z\right)\right)^2\right]$. 
		The estimators $(\widehat{\beta}, \widehat{K})$ turn out to hold a strong paralellism with linear regression least squares. Consider the observations \({Y} \in \mathbb{R}^n\) and the design matrix \({X} \in \mathbb{R}^{n \times p}\); then the least squares estimator $\hat{\beta}_{OL S}$ is any solution of the \textit{normal equations} ${X}^T {X}\hat{\beta}_{OL S}= {X}^T{Y} .$ Analogously, by setting both $p-$dimensional gradients of the objective function in (\ref{emprisk}) to zero we arrive at the next result. 
		
		\begin{Proposition}[GI normal equations]\label{normal}
			$(\widehat{\beta}, \widehat{K_{opt}})$ must satisfy 
			$$\begin{aligned}
				X^T X \widehat{\beta}+X^T\left(X-M_Z\right) \widehat{K_{opt}} & =X^T Y, \\
				X^T\left(X-M_Z\right) \widehat{\beta}+X^T\left(X-M_Z\right) \widehat{K_{opt}} & =\left(X-M_Z\right)^T Y. 
			\end{aligned}$$
		\end{Proposition}

		\noindent	By subtracting the second normal equation from the first one, we get $M^T_ZM_Z \widehat{\beta} = M_Z^TY$. If $M^T_ZM$ is invertible,  $\widehat{\beta}=(M^T_ZM_Z)^{-1}M_Z^TY$. By subtracting  $M^T_ZM_Z \widehat{\beta} = M_Z^TY$ to the second normal equation, we obtain $\widehat{K}_{opt}=(X^TX - M_Z^TM_Z)^{-1}X^T(Y-X\widehat{\beta})$. Using Remark \ref{remnz}, $\widehat{K}=\frac{1}{\sum_{z=1}^Z n_z}\left(X^TX - M^T_ZM_Z\right)\widehat{K_{opt}}=\frac{1}{	\sum_{z=1}^Z n_z}X^T(Y-X\widehat{\beta})$. These conclusions are assembled in the next theorem, of core importance in our development.

		\begin{theorem}[GI closed-form estimators] Assume that $M_Z^T M_Z$ is invertible. Then 	$$\widehat{\beta}=(M_Z^T M_Z )^{-1}M_Z^T Y \textrm{ and } 
		\widehat{K}  =\frac{1}{\sum_{z=1}^Z n_z}X^T(Y-X \widehat{\beta}). $$
			
		\end{theorem}
		
		\noindent Note that $\hat \beta$ resembles a two-stage least squares estimator that utilizes domain dummy indicators as instrumental variables \citep{long}. However, the key contribution of GI lies in incorporating $K_*$ for delivering optimal predictions adapted to a test domain, even when an asymptotically unbiased estimator of the causal parameter $\beta_*$ is employed. The existence and uniqueness of solutions to the normal equations is guaranteed under the same conditions for $\widehat{\mu}_1, \ldots, \widehat{\mu}_Z$ that Corollary \ref{multidentif} imposes on their population versions, that is $\operatorname{dim }\langle \widehat{\mu}_1, \ldots, \widehat{\mu}_Z\rangle = p.$

		% Let $X=(X_1,\ldots,X_n) \in \mathbb{R}^{n \times 1}$ , $Y=(Y_1,\ldots, Y_n) \in \mathbb{R}^{n \times 1}$  and $M= \left( \begin{array}{c|c|c}
			% 1_n   P_n X_1 &\cdots &1_n P_n X^p
			% \end{array}\right)$, where $1_n$ is the one of $\mathbb{R}^n$. Then

		\begin{example}
			Let $X=(X_1,\ldots,X_n) \in \mathbb{R}^{n }$ and $Y=(Y_1,\ldots, Y_n) \in \mathbb{R}^{n}$. $\widehat{\beta}$ and $\widehat{K}$ exist and are unique if invertibility of $M_1^TM_1$ is ensured. For this, it is enough to gather data from just one training environment such that the average of the covariate therein is not zero, agreeing with the univariate identifiability assumption in Theorem \ref{t1} as the aforementioned matrix is equal to $n\widehat{\mu}^2$. If our model is still 1D but includes an intercept, at least $2$ sources with different covariate averages must be used. The design matrix would be \texttt{cbind(1,X)}, where \texttt{X} is the $(n_1 + n_2) \times 1$ matrix concatenating the covariate observations $n_1,n_2$ of both sources. If the average of \(X\) is the same across both sources, then the population IV is not defined.
		\end{example}

		\section{Generative Invariance asymptotics}\label{sec:asy}
		
		The asymptotic properties of an estimator are crucial for making probabilistic statements about its behavior as the sample size grows large. Our proof strategies depart from standard theoretical arguments in the field. Typically, such reasoning chains rely on standard concentration inequalities for matrices with independent rows \citep{Vershynin_2012}. However, in the GI case, the rows involved in our matrices are not identically distributed. Specifically, the multivariate averages in each environment come from different distributions, which precludes the direct application of standard bounds.\footnote{Think of defining by columns $U:=\left(\widehat{\mu}_1 \mid \ldots \mid \widehat{\mu}_Z\right)$ and consider $A=U^T$ in Section 5.4.3 in \cite{Vershynin_2012}. The rows of $A$ (which are the columns of $U$, representing the multivariate averages in each environment) do not come from the same distribution.}
		Our first theorem establishes a concentration inequality, which serves to later prove asymptotic normality. Proofs for the results in this section are in the Supplement. We define $\Phi=\frac{1}{\sum_{z=1}^Z n_z} \sum_{z=1}^Zn_z \mu_z \mu_z^T$.

		\begin{theorem}[Finite-sample bound]\label{concineq}
			Consider multi-source data (\ref{data}), sampled i.i.d. from each source \(z\) such that each covariate observation is a \(p\)-dimensional random vector with mean \(\mu_z\) and covariance matrix \(\Sigma_z\), with sub-Gaussian components for all \(z\). Suppose the response noise \(\varepsilon_{Y_1},\ldots,\varepsilon_{Y_{n_z}}\) is also sub-Gaussian for all $z$. Let \(\lambda_z\) be \(\lambda_{\max}(\Sigma_z)\), \(\lambda_M\) be \(\max_{1, \ldots, Z} \lambda_z\), and \(\|\mu_M\|_2\) be \(\max_{1, \ldots, Z}\|\mu_z\|_2\). Assume that each source sample size $n_z > \frac{4}{\lambda_{\min }^2(\Phi)}\left((L^4+1) \lambda_M p+2(L^2 + 1)\left\|\mu_M\right\|_2 \sqrt{\lambda_M p}\right)^2$, where $L$ is the sub-Gaussian norm of a certain random variable not involving $p$ or $n$. Then, the following concentration inequality holds $$\mathbb{P}\left(\|\widehat{\beta}- \beta_{opt}\|_2 > t \right) \leq 12Z\exp{(-cp)},$$ where the slowest rate in $t$ involves $\max_z \sqrt \frac{p}{n_z}$ and $c$ is an absolute constant. 
		\end{theorem}

		\begin{Corollary}\label{cicausal}
			If the covariate means of the different sources span $\mathbb{R}^p$, then Theorem \ref{concineq} holds for $\beta_{opt}=\beta_*$, the causal parameter. 
		\end{Corollary}

		Next, we show that $\hat{\beta}$ is asymptotically normally distributed. Deriving the expression of the asymptotic covariance matrix has the practical consequence of enabling us to identify the necessary conditions for a direction to be in its null space. This is particularly useful for selecting and recombining data sources to maximize asymptotic efficiency when \(Z > p\). If \(Z = p\), all sources are required. Practical implications of knowing the asymptotic covariance matrix are discussed in Section \ref{real}, together with formal support in the Supplement. 
		
		\begin{theorem}[Asymptotic normality of $\hat{\beta}$]\label{asyn} Let the population and empirical mean vectors of the covariates corresponding to the different environments span $\mathbb{R}^p$. Suppose the random variables following the joint distribution of each source to be independent and have finite fourth order moments. Under the conditions in Theorem \ref{concineq},
			$$\left(\Phi^{-1} \sum_{z=1}^Z \frac{w_z^2}{n_z} \Omega_z \Phi^{-1}\right)^{-\frac{1}{2}}\left(\hat{\beta}-\beta_*\right) \quad \xrightarrow{d}  \quad \mathcal{N}\left(0, \operatorname{I}_p\right)$$
			\noindent with $\Omega_z:=\left(\mu_z^T \beta_*\right)^2 \Sigma_z-\mu_z^T \beta_*\left(K_* \mu_z^T+\mu_z K_*^T\right)+\mu_z \mu_z^T \sigma_Y^2.$
			
		\end{theorem}

		Let us examine the positive definiteness of the asymptotic covariance matrix of $\hat{\beta}$. The following result provides a starting point for an initial approximation of how to select the best data sources to feed our estimators, ensuring well-conditioning of the asymptotic covariance matrix of $\hat{\beta}$ when $Z > p$ environments are available. The practical implications of this result are further explored in the Supplement.
		
		\begin{Proposition} \label{omegabound}Let $z \in \{1,\ldots,Z\}$ and $\Omega_z$ be as in Theorem \ref{asyn}. Assume that causal Mahalanobis condition $z$ (\ref{schur_cond}) holds. Then for an arbitrary $0 \neq v \in \mathbb{R}^p$ 
			
			$$v^T \Omega_z v > \left(\mu_z^T \beta_* \sqrt{v^T \Sigma_z v}-\mu_z^T v \sigma\right)^2.$$
		\end{Proposition}
		
		\noindent Therefore, the causal Mahalanobis condition $z$ holding simultaneously for $z = 1, \ldots, Z$ formally confirms the positive definiteness of the asymptotic covariance matrix of $\beta$, as it is sufficient to verify this property for each $\Omega_z$, $z = 1, \ldots, Z$; providing a lower bound for its quadratic forms.

		We end this section comparing the asymptotic efficiency of GI with that of Causal Dantzig (CD) \citep{rothenhausler2019causal}, when $p=1$ and $Z=2$. The asymptotic variance of the CD estimator of $\beta_*$ is $V_{CD} = \frac{1}{n} \frac{2}{\mathbb{E}X_1^2 - \mathbb{E}X_2^2} \left( \operatorname{Var}(X_1 \varepsilon_{Y,1}) + \operatorname{Var}(X_2 \varepsilon_{Y,2}) \right)$, whereas that of the GI estimator of $\beta_*$ is
		$V_{GI} = \frac{2}{n}
		\frac{(\mu_1^2 + \mu_2^2)(\sigma_Y^2 - 2\beta_* K_*) + \beta_*^2 \left(\mu_1^2 \operatorname{Var}(X_1) + \mu_2^2 \operatorname{Var}(X_2) \right)}{(\mu_1^2 + \mu_2^2)^2}$. 
		Even though an additive intervention on $X_1$ being non-zero and having a positive second-order moment would suffice for identifiability of CD, according to Theorem 2 in \cite{rothenhausler2019causal}, it could be that the denominator of $V_{CD}$ is close to zero, inflating the whole fraction. Assume the mean of the intervention on $X$ across environments is zero (this still preserves identifiability of CD).  Therefore, a small denominator of $V_{CD}$ is equivalent to $\operatorname{Var}(X_1) \approx \operatorname{Var}(X_2)$ because we have $\mu_1=\mu_2=:\mu$. In this setting, $V_{GI} \approx \frac{1}{n} \frac{\sigma_Y^2 - 2\beta_* K_* + \Var(X)\beta^2_*}{ \mu^2}$. 	Thus, the asymptotic efficiency of GI does not suffer when $\operatorname{Var}(X_1) \approx \operatorname{Var}(X_2)$, as long as $\mu$  is large enough in absolute value. Moreover, causal Mahalanobis conditions 1 and 2 would mean here that $K^2_* < \sigma^2_Y\Var X$, ensuring $V_{GI} > \frac{1}{n} \frac{\left(\beta_* \sqrt{\Var (X)} -\frac{ K_*}  {\sqrt{\Var (X) }}\right)^2}{\mu^2} \geq 0$.

		\section{Numerical examples}
		
		This section empirically demonstrates the superior performance of GI. We begin with an illustrative univariate medical example where data is sourced from different hospitals, highlighting the limitations of traditional methods and the improvements provided by GI. Following this, we perform an empirical comparison with DRIG, a method based on finite robustness. Next, we benchmark GI against well-known domain adaptation methods, including DIP and CIRM, to emphasize its superior target risk performance even in misspecified settings. Subsequently, the comparison is extended to the IV and CD approaches in setups favoring such methods, showcasing GI’s robustness and predictive accuracy. An analysis of GI’s distributional optimality under gaussianity is also provided, using energy distance metrics (see Supplement) to quantify its ability to emulate test distributions across interventions. Finally, we cover an application of GI to real-world data from the SPRINT trial, presenting its advantages in clinical decision support. 
		\subsection{The problem of unsupervised domain adaptation: an example}
		The estimated glomerular filtration rate (EGFR) measures how well kidneys filter waste from blood. Systolic blood pressure (SBP) is the pressure exerted on artery walls when the heart beats. We want to develop a predictive model of EGFR using SBP. Hypertension (high SBP) is a known risk factor for chronic kidney disease (CKD), potentially leading to a decline in kidney function over time. By accurately predicting EGFR from SBP, healthcare providers could identify individuals at risk of CKD earlier and implement preventive measures or treatments to mitigate future kidney damage, thereby improving patient outcomes and reducing the burden of renal disease on the healthcare system.
		
		The solid black line fit in Figure \ref{in1} (left) is solely based on data from hospital 48 plotted as solid black dots. The \texttt{lm} summary in \texttt{R} outputs a \textit{p}-value of $0.06$ for the slope estimator. Subsequently, we realize there is access to data from hospital 22 (white squares in Figure \ref{in1}). We fit a second linear model, represented by the solid red line in Figure \ref{in1} (left), to the pooled data from both hospitals, 48 and 22. This fit aims to predict SBP values from a third hospital, hospital 15, represented as blue marks on the horizontal axis. Predictions given by GI for hospital 15 are shown in blue (GI did not have access to the response values from test hospital 15). Despite not using the triangles (only the rug on the horizontal axis), GI is able to produce much better predictions for the third hospital than ordinary least squares.

		\subsection{Empirical comparison with finite robustness}
		In order to evaluate the efficacy of GI, we conduct a series of simulations comparing it with the state-of-the-art method known as {Distributional Robustness via Invariant Gradients} (DRIG) \citep{shen2023causality}. It is a regularized empirical risk minimizer based on a gradient invariance condition across environments, offering robustness against interventions through a controlled regularization parameter $\gamma$. For instance, DRIG corresponds to performing OLS on a pre-selected environment when $\gamma = 0$, which is represented by the red dots in Figure \ref{firstfig}. The robustness of DRIG, particularly against a variety of intervention strengths and directions, sets a high benchmark for new methodologies. We use a simulation setup that is natural for DRIG, and even counterproductive for GI as it does not guarantee that the average vectors of covariates across environments span $\mathbb{R}^p$. Despite this, we showcase that GI surpasses DRIG in terms of prediction mean squared error in this example. The results, displayed in Figure \ref{simdrig}, show that even without enforcing the assumptions necessary for our estimator's optimal performance, our approach uniformly outperforms DRIG across different perturbation strengths and \(\gamma\) values in this example. We refer the reader to the Supplement for technical details about the simulation setup.

		{\subsection{Empirical comparison with unsupervised domain adaptation methods}}
		
		We conduct a numerical comparison within a misspecified setting for GI, benchmarking it against methods from the unsupervised domain adaptation literature. Interestingly, GI consistently outperforms them even when Assumption (\ref{firstscm}) is violated. The domain adaptation methods used here are as follows: \texttt{SrcPool}, which pools all training sources together and performs OLS; \texttt{DIPweigh}, which aligns the source and target empirical distributions in a low-dimensional latent space, making it suitable for multi-source data; \texttt{CIP}, which leverages label information from multiple sources; \texttt{CIRMweigh}, which combines \texttt{DIP} and \texttt{CIP}; and \texttt{CIRM<>weigh}, an oracle that has knowledge of the true causal covariate indices. For a deeper review on the domain adaptation methods used here, please refer to the Supplement.

		We base our comparison on target risk over 1000 random coefficient generations. The experimental setup mirrors simulation vii) from Section 5.2 in \cite{chen}, with $Z=11$ training environments (to account for the intercept, given a dimensionality of $p = 10$, requiring one additional source environment for GI identification). In the training source $z$, each data point is generated i.i.d. according to the SCM
		$
		\begin{bmatrix}
			X_z \\ 
			Y_z
		\end{bmatrix}
		=
		\begin{bmatrix}
			\mathbf{B} & b \\ 
			\omega^\top & 0
		\end{bmatrix}
		\begin{bmatrix}
			X_z \\ 
			Y_z
		\end{bmatrix}
		+ \text{additive intervention}_z + \text{endogenous noise}_z.
		$
		Here, $b = 0$ corresponds to a \textit{causal} setting (arrows directed toward $Y$, coherent with the SCM assumed for GI), $\omega = 0$ represents an \textit{anticausal} setting (arrows directed outward from $Y$), and other configurations are categorized as \textit{mixed-causal-anticausal}. In our simulations, we use a mixed-causal-anticausal setting with $b^j = 0$ for odd $j$ and $\omega^j = 0$ for even $j$. This setup is intentionally misspecified for GI, as the assumed SCM (\ref{firstscm}) considers arrows directed inward to $Y$. We intervene on the response by altering its variance as $\sigma_{Y,z}= \sigma_Y \mathcal{U}(1,2)$. Remark \ref{diff_var} clarifies why identifiability of GI is still possible under this class of interventions on $Y$. Endogenous noise accounting for hidden confounding is generated by (randomly) linearly combining the noise components of $X$ and adding exogenous noise scaled so that $\sigma_{Y,z}$ is respected. Results are visible in Figure \ref{fig:diag}\footnote{Figure \ref{fig:diag} was produced using code adapted from  \texttt{https://github.com/yuachen/CausalDA}.}.

		%\begin{table}[h!]
		%	\centering
		%	\begin{tabular}{lccccc}
			%		\toprule
			%		Method Name           & \texttt{OLS} & \texttt{DIPweigh} &\texttt{CIP}& \texttt{CIRMweigh} & Oracle \texttt{CIRMweigh} \\
			%		\midrule
			%		GI better (out of 100 runs) & 84  & 93           & 84   & 75   &76       \\
			%		\bottomrule
			%	\end{tabular}
		%	\caption{Comparison of methods based on the number of times GI outperforms each method over 100 random coefficient \textbf{B} data generations. \texttt{OLS} stands for pooling all training sources together and then performing OLS. Oracle \texttt{CIRMweigh} knows the indexes of the true causal covariates. }
		%	\label{tab:gi_comparison}
		%\end{table}

		{	\subsection{Empirical comparison with IV and Causal Dantzig}}
		
		We exactly replicate the simulation setup described in model (25) of \cite{rothenhausler2019causal}, where the authors compare CD with IV in a scenario where the environment-conditional mean of the covariate ($p=1$)  changes only slightly. This leads to technical challenges for the IV approach. However, as the second order moment of the covariate changes substantially across environments, CD is well-defined at the population level. This results in favorable asymptotic properties for estimating $\beta_*$, but does not necessarily determine the quality of predictions in an unseen environment. Interestingly, even a standard IV approach performs better than CD in this context. In turn, GI consistently outperforms both approaches across all sample sizes, as Table \ref{simdantz} shows.

		\subsection{Distributional optimality of Generative Invariance}

		Corollary \ref{copy} asserts that, under gaussianity, the generator $g_{K_*}(\cdot,\cdot)$ has the same distribution as the test environment endogenous residuals. To evaluate this, we create a simulation setup that allows us to compare the performance of GI against just adding exogenously sampled noise for each fixed $x \in \mathbb{R}$ to $x\widehat{\beta}_{OLS}$ (``OLS'') or even to $x\widehat{\beta}$ (``Causal parameter'') in Figure \ref{fig:energy}. The comparison is based on how well these methods replicate the joint distribution of several test environments. We use the energy distance \citep{szekely2017energy} as a metric.

		Results are visible in Figure \ref{fig:energy}, while a detailed description of the setup can be found in the Supplement. Notice how the causal parameter performs worse than OLS when the test covariate distribution is similar to that of the training set. This happens at $\mathbb{E}V^2 = 2$ because the training covariate distribution had been perturbed by a Gaussian random variable with mean and variance equal to 1 and, therefore,	 second order moment equal to 2. Moreover, observe how OLS degrades arbitrarily as the strength of interventions grows- there is a log scale in the horizontal axis. GI performs uniformly well across interventions.

		\subsection{Application: the SPRINT trial}\label{real}
		
		This section performs an acid test to GI checking its performance in a fresh sample constituted by data sources not seen during training so as to showcase its superiority in terms of test MSE in comparison to the state-of-the-art \citep{shen2023causality}, no matter what the value of its hyperparameter is, and OLS. We will see as well how the consequences of our asymptotic normality result (Theorem \ref{asyn}) lead to practical tools that serve to automatically decide the sources used to train GI optimally. 
		
		%	The NIH’s Systolic Blood Pressure Intervention Trial (SPRINT) was conducted to inform the new blood pressure medication guidelines in the US by testing the effects that a lower blood pressure target has on reducing heart disease risk. Observational studies had shown that individuals with lower SBP levels had fewer complications and deaths due to cardiovascular disease (CVD). Building on this observation, SPRINT was designed to test the effects of a lower blood pressure target on reducing heart disease risk. Specifically, SPRINT aimed to compare treating high blood pressure to a target SBP goal of less than 120 mmHg against treating to a goal of less than 140 mmHg.
		
		In the context of the SPRINT trial \citep{sprint}, the variable  RISK10YRS  represents the 10-year risk of cardiovascular events, calculated using the Framingham Risk Score \citep{wilson1998prediction}, a widely used metric to quantify the 10-year cardiovascular risk. It is crucial as it encapsulates the long-term cardiovascular risk profile of participants, providing a measurable outcome to assess the effectiveness of different blood pressure targets tested in the trial. Understanding the factors influencing RISK10YRS  is essential for developing predictive models that can guide clinical decisions and policy-making. In particular, building a predictive model for it as a function of  AGE, SBP and BMI is of major interest as these variables are well-established  risk factors for cardiovascular disease \citep{ordovas}. 
		
		However, the potential impact of omitted confounders on the predictive model cannot be understated. Confounders such as smoking status, aspirin use, sex and diabetes can influence both the independent and dependent variables. For instance, smoking and aspirin use are known to affect cardiovascular risk independently of blood pressure. The omission of these confounders from the model leads to biased estimates of the effects of AGE, SBP and BMI on RISK10YRS, undermining the validity of the predictive model. This is where GI comes into play: the parameter $K \in \mathbb{R}^3$ takes care of accounting for noise endogeneity due to excluding variables related to the covariates yet simultaneously influencing RISK10YRS. Since there are three covariates and an intercept in this problem, according to Corollary \ref{multidentif} we need data from at least $p=4$ different sources and that the multivariate average vectors of the covariates are linearly independent.
		
		SPRINT was conducted at 102 clinical sites in the US. We estimate the energy distance matrix of the covariates across sites and average it in the direction of rows (columns would give the same output because of symmetry). Ranking these averages highlights the sites that are distributionally the most distant with respect to the others in terms of patient characteristics. We choose the 10 most peculiar sites according to this criterion for testing purposes. For training, we consider the sites with more than 100 observations after excluding the test ones, 19 in total. We train OLS and DRIG (for several values of its hyperparameter $\gamma$) using the data from the 19 training environments. Strictly speaking, we do not use all these sources for training GI, but look for the top combinations of $p+1=4$ sources, producing as many GI estimates as sets of sources we choose and then average across them the predictions given by our estimator. We start by computing $4$-dimensional covariate average vectors within each site. Then choose a subset of 4 sources among the 19 environments, stack within a 4x4 matrix the columns corresponding to the averages and compute its determinant, repeating this procedure for the $\binom{19}{4}=3876$ possible combinations (0.037 seconds in \texttt{R}) and sort the combinations in decreasing determinant orders. We consider the best $B=100$ combinations, obtain $B$ GI estimators and exploit them as in the Supplement. 
		
		For each test environment, $\frac{B}{2}$ predictions are launched per covariate observation and then average them. We estimate OLS pooling the 19 training environments and DRIG for $\gamma \in \{0.05,0.1,0.3,0.5\}$, with the median MSE computed across the 10 test environments. Results are displayed in Table \ref{numres} and Figure \ref{testplots}, the latter offering a picture of predictive behavior in some of the test sites. The most favorable value for DRIG $\gamma=0.1$ was picked.

		Observe that we are not putting our estimator in a privileged position in terms of the amount of training data: the six approaches use the same training environments. A different story is that GI automatically discards some of them and combines the remaining ones following a procedure that maximizes efficiency (Theorem \ref{asyn}).
		
		For a fixed test environment, the pre-aggregated predictions given by GI are not carried out through independent draws of the underlying estimators, since different source combinations do not necessarily have an empty intersection. However, we believe that the reason this procedure still outperforms the baseline algorithms is related to phenomenology similar to that described for the bootstrap: as long as the base estimator has a CLT, the resampling distribution is well-behaved even though it is built upon dependent samples \citep{kosorok}.

		One advantage of GI compared to DRIG is that the latter has a strong dependence on the so-called \textit{reference environment}, as their estimator barely considers data from other sources when \(\gamma \approx 0\). On top of outperforming state-of-the-art estimators, ours is asymptotically unbiased: we can target the causal parameter and still use it to provide optimal predictions instead of committing identifiability bias.

		\section{Conclusions}
		
		We have presented a novel approach to domain adaptation under hidden confounding that challenges the prevailing notion that, in some situations, identifiability must be compromised for robust prediction. Our proposal, Generative Invariance, represents a novel perspective in statistical causal inference, since it can be interpreted as a \textit{causal} reconciliation of \cite{breiman}’s “two cultures”, bridging the gap between identifiability and prediction. In addition, it offers significant advantages including less performance degradation as strength of test perturbations grows, not needing multi-source data in the one-dimensional case, and the elimination of hyperparameter tuning.  We make no assumptions on how much the training and test distributions overlap, as often needed in the literature, see e.g. \cite{garciameixide}.

		Recent work in the causal representation learning literature still requires the assumption that the number of environments is at least equal to the number of latent nodes \citep{syrcrl}. Recently, \cite{aragam} enabled causal representation learning with fewer environments by identifying latent nodes with shifted mechanisms rather than reconstructing the full causal structure, comparing noise components across environments to detect shifts. Weak identification in nonparametric IV is an actively ongoing research stream \citep{newey}}. 
	
	It would be conceivable to model the population parameters $\mu_1, \ldots, \mu_Z$ as being drawn from a common prior distribution. Thus, we would be effectively performing an empirical Bayes analysis estimating the second-order moment of the prior (precisely equal to $\sum_{z=1}^Z \widehat{\mu}_z \widehat{\mu}_z^T= M^T_ZM_Z$ if prior is zero mean). The hypothetical \textit{posterior} distribution of $\beta_*$ would therefore have cross-world \citep{richardson2013single} implications and would endow summation over $z$ with an asymptotic meaning.

	Promising avenues for further research include extending our approach to high-dimensional settings, exploring nonlinear models and obtaining approximate identifiability results by bounding $\lVert \beta_{opt} - \beta_* \rVert_{\infty}$ and $\lVert K_{opt} - K_* \rVert_{\infty}$ when there are interventions in the latent variables. Valid inference after source combination selection, as discussed in the Supplement based on the asymptotic confidence intervals for $\beta_*$ provided by Theorem \ref{asyn} opens an interesting research line, as well. Also, the class $\mathcal{M}$ looks like a semiparametric model when replacing $\beta^Tx$ by a richer family of functions. The part of the model concerning $g_K$ is always going to be parametric, having the dimension as the covariates. These exciting open questions not only promise to keep advancing our current understanding of causal inference but also to enhance practical methodologies employed in the field of statistics.

	\bibliographystyle{biorefs}
	%\bibliography{Bibliography-MM-MC}

\begin{thebibliography}{99}

\bibitem[Angrist \emph{and others}(1996)Angrist, Imbens and
  Rubin]{angrist1996identification}
\textsc{Angrist, Joshua~D, Imbens, Guido~W and Rubin, Donald~B}. (1996).
\newblock Identification of causal effects using instrumental variables.
\newblock {\em Journal of the American Statistical
  Association\/}~\textbf{91}(434), 444--455.

\bibitem[Arjovsky \emph{and others}(2020)Arjovsky, Bottou, Gulrajani and
  Lopez-Paz]{Arjovsky2019}
\textsc{Arjovsky, Martin, Bottou, Léon, Gulrajani, Ishaan and Lopez-Paz,
  David}. (2020).
\newblock Invariant risk minimization.

\bibitem[Baktashmotlagh \emph{and others}(2013)Baktashmotlagh, Harandi, Lovell
  and Salzmann]{dip13}
\textsc{Baktashmotlagh, Mahsa, Harandi, Mehrtash~T., Lovell, Brian~C. and
  Salzmann, Mathieu}. (2013).
\newblock Unsupervised domain adaptation by domain invariant projection.
\newblock In:  {\em 2013 IEEE ICCV\/}. pp.\  769--776.

\bibitem[Bennett \emph{and others}(2023)Bennett, Kallus, Mao, Newey, Syrgkanis
  and Uehara]{newey}
\textsc{Bennett, Andrew, Kallus, Nathan, Mao, Xiaojie, Newey, Whitney,
  Syrgkanis, Vasilis and Uehara, Masatoshi}. (2023).
\newblock Inference on strongly identified functionals of weakly identified
  functions.

\bibitem[Breiman(2001)Breiman]{breiman}
\textsc{Breiman, Leo}. (2001).
\newblock {Statistical Modeling: The Two Cultures (with comments and a
  rejoinder by the author)}.
\newblock {\em Statistical Science\/}~\textbf{16}(3), 199 -- 231.

\bibitem[Cauchois \emph{and others}(2024)Cauchois, Gupta, Ali and Duchi]{duchi}
\textsc{Cauchois, Maxime, Gupta, Suyash, Ali, Alnur and Duchi, John~C.} (2024).
\newblock Robust validation: Confident predictions even when distributions
  shift.
\newblock {\em Journal of the American Statistical
  Association\/}~\textbf{0}(0), 1--66.

\bibitem[Chen \emph{and others}(2024)Chen, Bello, Locatello, Aragam and
  Ravikumar]{aragam}
\textsc{Chen, Tianyu, Bello, Kevin, Locatello, Francesco, Aragam, Bryon and
  Ravikumar, Pradeep}. (2024).
\newblock Identifying general mechanism shifts in linear causal
  representations.

\bibitem[Chen and Buehlmann(2021)Chen and Buehlmann]{chen}
\textsc{Chen, Yuansi and Buehlmann, Peter}. (2021).
\newblock Domain adaptation under structural causal models.
\newblock {\em Journal of Machine Learning Research\/}~\textbf{22}(261), 1--80.

\bibitem[Chernozhukov \emph{and others}(2018)Chernozhukov, Chetverikov,
  Demirer, Duflo, Hansen, Newey and Robins]{chernodoubly}
\textsc{Chernozhukov, Victor, Chetverikov, Denis, Demirer, Mert, Duflo, Esther,
  Hansen, Christian, Newey, Whitney and Robins, James}. (2018, 01).
\newblock {Double/debiased machine learning for treatment and structural
  parameters}.
\newblock {\em The Econometrics Journal\/}~\textbf{21}(1), C1--C68.

\bibitem[Christiansen \emph{and others}(2021)Christiansen, Pfister, Jakobsen,
  Gnecco and Peters]{christiansen2021causal}
\textsc{Christiansen, Rune, Pfister, Niklas, Jakobsen, Martin~Emil, Gnecco,
  Nicola and Peters, Jonas}. (2021).
\newblock A causal framework for distribution generalization.
\newblock {\em IEEE Transactions on Pattern Analysis and Machine
  Intelligence\/}~\textbf{44}(10), 6614--6630.

\bibitem[Didelez \emph{and others}(2006)Didelez, Dawid and
  Geneletti]{didelez2006direct}
\textsc{Didelez, Vanessa, Dawid, A.~Philip and Geneletti, Sara}. (2006).
\newblock Direct and indirect effects of sequential treatments.
\newblock In:  {\em Proceedings of the 22nd Annual Conference on Uncertainty in
  Artificial Intelligence\/}. Corvallis. pp.\  138--146.

\bibitem[Foster and Syrgkanis(2023)Foster and Syrgkanis]{syrg}
\textsc{Foster, Dylan~J. and Syrgkanis, Vasilis}. (2023).
\newblock {Orthogonal statistical learning}.
\newblock {\em The Annals of Statistics\/}~\textbf{51}(3), 879 -- 908.

\bibitem[Goodfellow \emph{and others}(2014)Goodfellow, Pouget-Abadie, Mirza,
  Xu, Warde-Farley, Ozair, Courville and Bengio]{goodfellow2014generative}
\textsc{Goodfellow, Ian, Pouget-Abadie, Jean, Mirza, Mehdi, Xu, Bing,
  Warde-Farley, David, Ozair, Sherjil, Courville, Aaron and Bengio, Yoshua}.
  (2014).
\newblock Generative adversarial nets.
\newblock {\em Advances in neural information processing
  systems\/}~\textbf{27}.

\bibitem[Haavelmo(1943)Haavelmo]{Haavelmo1943}
\textsc{Haavelmo, T.} (1943).
\newblock The statistical implications of a system of simultaneous equations.
\newblock {\em Econometrica\/}~\textbf{11}, 1--12.

\bibitem[Hall(2005)Hall]{hall2005generalized}
\textsc{Hall, Alastair~R.} (2005).
\newblock {\em Generalized method of moments\/}. Oxford University Press.

\bibitem[He and Geng(2008)He and Geng]{he2008active}
\textsc{He, Yang-Bo and Geng, Zhi}. (2008).
\newblock Active learning of causal networks with intervention experiments and
  optimal designs.
\newblock {\em Journal of Machine Learning Research\/}~\textbf{9}, 2523--2547.

\bibitem[Henzi \emph{and others}(2024)Henzi, Shen, Law and Bühlmann]{henzinva}
\textsc{Henzi, Alexander, Shen, Xinwei, Law, Michael and Bühlmann, Peter}.
  (2024, 11).
\newblock Invariant probabilistic prediction.
\newblock {\em Biometrika\/}, asae063.

\bibitem[Horn and Johnson(2012)Horn and Johnson]{john}
\textsc{Horn, Roger~A. and Johnson, Charles~R.} (2012).
\newblock {\em Matrix Analysis\/}, 2 edition. Cambridge University Press.

\bibitem[Javanmard and Montanari(2014)Javanmard and Montanari]{montanari}
\textsc{Javanmard, Adel and Montanari, Andrea}. (2014).
\newblock Confidence intervals and hypothesis testing for high-dimensional
  regression.
\newblock {\em Journal of Machine Learning Research\/}~\textbf{15}(82),
  2869--2909.

\bibitem[Jeong and Rothenhäusler(2024)Jeong and
  Rothenhäusler]{jeong2024outofdistribution}
\textsc{Jeong, Yujin and Rothenhäusler, Dominik}. (2024).
\newblock Out-of-distribution generalization under random, dense distributional
  shifts.

\bibitem[Jin and Syrgkanis(2024)Jin and Syrgkanis]{syrcrl}
\textsc{Jin, Jikai and Syrgkanis, Vasilis}. (2024).
\newblock Learning causal representations from general environments:
  Identifiability and intrinsic ambiguity.

\bibitem[Kosorok(2008)Kosorok]{kosorok}
\textsc{Kosorok, Michael~R}. (2008).
\newblock {\em Introduction to empirical processes and semiparametric
  inference\/}, Volume~61. Springer.

\bibitem[Kostin \emph{and others}(2024)Kostin, Gnecco and Yang]{kostin}
\textsc{Kostin, Julia, Gnecco, Nicola and Yang, Fanny}. (2024).
\newblock Achievable distributional robustness when the robust risk is only
  partially identified.
\newblock In:  {\em The Thirty-eighth Annual Conference on Neural Information
  Processing Systems\/}.

\bibitem[Krueger \emph{and others}(2021)Krueger, Caballero, Jacobsen, Zhang,
  Binas, Zhang, Le~Priol and Courville]{krueger2021out}
\textsc{Krueger, David, Caballero, Ethan, Jacobsen, Joern-Henrik, Zhang, Amy,
  Binas, Jonathan, Zhang, Dinghuai, Le~Priol, Remi and Courville, Aaron}.
  (2021).
\newblock Out-of-distribution generalization via risk extrapolation (rex).
\newblock In:  {\em ICML\/}. PMLR. pp.\  5815--5826.

\bibitem[Lin \emph{and others}(2022)Lin, Dong, Wang and Zhang]{lin2022bayesian}
\textsc{Lin, Yong, Dong, Hanze, Wang, Hao and Zhang, Tong}. (2022).
\newblock Bayesian invariant risk minimization.
\newblock In:  {\em Proceedings of the IEEE/CVF Conference on Computer Vision
  and Pattern Recognition\/}. pp.\  16021--16030.

\bibitem[Long \emph{and others}(2023)Long, Zhu, Do and Ha]{long}
\textsc{Long, James~P., Zhu, Hongxu, Do, Kim-Anh and Ha, Min~Jin}. (2023).
\newblock {Estimating causal effects with hidden confounding using instrumental
  variables and environments}.
\newblock {\em Electronic Journal of Statistics\/}~\textbf{17}(2), 2849 --
  2879.

\bibitem[Magliacane \emph{and others}(2017)Magliacane, van Ommen, Claassen,
  Bongers, Versteeg and Mooij]{Magliacane2017}
\textsc{Magliacane, S., van Ommen, T., Claassen, T., Bongers, S., Versteeg, P.
  and Mooij, J.~M.} (2017).
\newblock Domain adaptation by using causal inference to predict invariant
  conditional distributions.
\newblock In:  {\em Neural Information Processing Systems\/}.

\bibitem[Meixide and Matabuena(2025)Meixide and Matabuena]{garciameixide}
\textsc{Meixide, Carlos~García and Matabuena, Marcos}. (2025).
\newblock Causal survival embeddings: Non-parametric counterfactual inference
  under right-censoring.
\newblock {\em Statistical Methods in Medical Research\/}~\textbf{34}(3),
  574--593.
\newblock PMID: 39930905.

\bibitem[Ordovas \emph{and others}(2023)Ordovas, Rios-Insua, Santos-Lozano,
  Lucia, Torres, Kosgodagan and Camacho]{ordovas}
\textsc{Ordovas, J.M., Rios-Insua, D., Santos-Lozano, A., Lucia, A., Torres,
  A., Kosgodagan, A. and Camacho, J.M.} (2023).
\newblock A {Bayesian} network model for predicting cardiovascular risk.
\newblock {\em Computer Methods and Programs in Biomedicine\/}~\textbf{231},
  107405.

\bibitem[Pan \emph{and others}(2011)Pan, Tsang, Kwok and Yang]{tca}
\textsc{Pan, Sinno~Jialin, Tsang, Ivor~W., Kwok, James~T. and Yang, Qiang}.
  (2011).
\newblock Domain adaptation via transfer component analysis.
\newblock {\em IEEE Transactions on Neural Networks\/}~\textbf{22}(2),
  199--210.

\bibitem[Pearl(2009)Pearl]{Pearl_2009}
\textsc{Pearl, Judea}. (2009).
\newblock {\em Causality\/}, 2 edition. Cambridge University Press.

\bibitem[Peters \emph{and others}(2016)Peters, B{\"u}hlmann and
  Meinshausen]{peters2016causal}
\textsc{Peters, Jonas, B{\"u}hlmann, Peter and Meinshausen, Nicolai}. (2016).
\newblock Causal inference by using invariant prediction: identification and
  confidence intervals.
\newblock {\em J. of the Royal Statistical Society Series B: Statistical
  Methodology\/}~\textbf{78}(5), 947--1012.

\bibitem[Peters \emph{and others}(2017)Peters, Janzing and
  Sch{\"o}lkopf]{peters2017elements}
\textsc{Peters, Jonas, Janzing, Dominik and Sch{\"o}lkopf, Bernhard}. (2017).
\newblock {\em Elements of causal inference: foundations and learning
  algorithms\/}. The MIT Press.

\bibitem[Richardson and Robins(2013)Richardson and
  Robins]{richardson2013single}
\textsc{Richardson, Thomas and Robins, James~M.} (2013).
\newblock Single world intervention graphs (swigs): A unification of the
  counterfactual and graphical approaches to causality.
\newblock {\em Working Paper} 128, University of Washington.

\bibitem[Robins(1986)Robins]{robins86}
\textsc{Robins, James~M.} (1986).
\newblock A new approach to causal inference in mortality studies with a
  sustained exposure period—application to control of the healthy worker
  survivor effect.
\newblock {\em Mathematical Modelling\/}~\textbf{7}, 1393--1512.

\bibitem[Rothenh{\"a}usler \emph{and others}(2019)Rothenh{\"a}usler,
  B{\"u}hlmann and Meinshausen]{rothenhausler2019causal}
\textsc{Rothenh{\"a}usler, Dominik, B{\"u}hlmann, Peter and Meinshausen,
  Nicolai}. (2019).
\newblock Causal dantzig: Fast inference in linear structural equation models
  with hidden variables under additive interventions.
\newblock {\em The Annals of Statistics\/}~\textbf{47}(3), 1688--1722.

\bibitem[Rothenh{\"a}usler \emph{and others}(2021)Rothenh{\"a}usler,
  Meinshausen, B{\"u}hlmann and Peters]{rothenhausler2021anchor}
\textsc{Rothenh{\"a}usler, Dominik, Meinshausen, Nicolai, B{\"u}hlmann, Peter
  and Peters, Jonas}. (2021).
\newblock Anchor regression: Heterogeneous data meet causality.
\newblock {\em Journal of the Royal Statistical Society Series B: Statistical
  Methodology\/}~\textbf{83}(2), 215--246.

\bibitem[Rubin(1974)Rubin]{rubin74}
\textsc{Rubin, Donald~B.} (1974).
\newblock Estimating causal effects of treatments in randomized and
  nonrandomized studies.
\newblock {\em Journal of Educational Psychology\/}~\textbf{66}(5), 688--701.

\bibitem[Saengkyongam \emph{and others}(2022)Saengkyongam, Henckel, Pfister and
  Peters]{saengkyongam2022exploiting}
\textsc{Saengkyongam, Sorawit, Henckel, Leonard, Pfister, Niklas and Peters,
  Jonas}. (2022).
\newblock Exploiting independent instruments: Identification and distribution
  generalization.
\newblock In:  {\em International Conference on Machine Learning\/}. PMLR. pp.\
   18935--18958.

\bibitem[Shen \emph{and others}(2023)Shen, B{\"u}hlmann and
  Taeb]{shen2023causality}
\textsc{Shen, Xinwei, B{\"u}hlmann, Peter and Taeb, Armeen}. (2023).
\newblock Causality-oriented robustness: exploiting general additive
  interventions.

\bibitem[Sokol \emph{and others}(2014)Sokol, Maathuis and Falkeborg]{ica}
\textsc{Sokol, Alexander, Maathuis, Marloes~H. and Falkeborg, Benjamin}.
  (2014).
\newblock {Quantifying identifiability in independent component analysis}.
\newblock {\em Electronic Journal of Statistics\/}~\textbf{8}(1), 1438 -- 1459.

\bibitem[Splawa-Neyman(1990)Splawa-Neyman]{splawa}
\textsc{Splawa-Neyman, J.} (1990).
\newblock On the application of probability theory to agricultural experiments.
  essay on principles. section 9.
\newblock {\em Statistical Science\/}~\textbf{5}(4), 465--472.
\newblock Translated from the Polish and edited by D. M. Dabrowska and T. P.
  Speed.

\bibitem[Sz{\'e}kely and Rizzo(2017)Sz{\'e}kely and Rizzo]{szekely2017energy}
\textsc{Sz{\'e}kely, G{\'a}bor~J and Rizzo, Maria~L}. (2017).
\newblock The energy of data.
\newblock {\em Annual Review of Statistics and Its Application\/}~\textbf{4},
  447--479.

\bibitem[The SPRINT Research~Group(2015)The SPRINT Research~Group]{sprint}
\textsc{The SPRINT Research~Group, National Insitute~Health}. (2015).
\newblock A randomized trial of intensive versus standard blood-pressure
  control.
\newblock {\em New England Journal of Medicine\/}~\textbf{373}(22), 2103--2116.

\bibitem[van~de Geer(2010)van~de Geer]{vandegeer2010empirical}
\textsc{van~de Geer, Sara~A.} (2010, January).
\newblock {\em Empirical Processes in M-Estimation\/}, Cambridge Series in
  Statistical and Probabilistic Mathematics. Cambridge, UK: Cambridge
  University Press.
\newblock Rijksuniversiteit Leiden, The Netherlands.

\bibitem[van~der Laan and Rubin(2006)van~der Laan and Rubin]{laan06}
\textsc{van~der Laan, Mark~J. and Rubin, Daniel}. (2006).
\newblock Targeted maximum likelihood learning.
\newblock {\em The International Journal of Biostatistics\/}~\textbf{2}, 1--40.

\bibitem[Vansteelandt and Morzywołek(2023)Vansteelandt and
  Morzywołek]{vansteelandt2023orthogonal}
\textsc{Vansteelandt, Stijn and Morzywołek, Paweł}. (2023).
\newblock Orthogonal prediction of counterfactual outcomes.

\bibitem[Vershynin(2012)Vershynin]{Vershynin_2012}
\textsc{Vershynin, Roman}. (2012).
\newblock {\em Introduction to the non-asymptotic analysis of random
  matrices\/}. Cambridge University Press, p. 210–268.

\bibitem[Vershynin(2018)Vershynin]{verhdp}
\textsc{Vershynin, Roman}. (2018).
\newblock {\em High-Dimensional Probability: An Introduction with Applications
  in Data Science\/}, Cambridge Series in Statistical and Probabilistic
  Mathematics. Cambridge University Press.

\bibitem[Wilson \emph{and others}(1998)Wilson, D'Agostino, Levy, Belanger,
  Silbershatz and Kannel]{wilson1998prediction}
\textsc{Wilson, Peter~WF, D'Agostino, Ralph~B, Levy, Daniel, Belanger,
  Albert~M, Silbershatz, H and Kannel, William~B}. (1998).
\newblock Prediction of coronary heart disease using risk factor categories.
\newblock {\em Circulation\/}~\textbf{97}(18), 1837--1847.

\bibitem[Wright(1928)Wright]{wright1928tariff}
\textsc{Wright, Philip~Green}. (1928).
\newblock {\em The tariff on animal and vegetable oils\/}, Number~26.
  Macmillan.

\bibitem[Zhang \emph{and others}(2013)Zhang, Sch{\"o}lkopf, Muandet and
  Wang]{Zhang2013}
\textsc{Zhang, K., Sch{\"o}lkopf, B., Muandet, K. and Wang, Z.} (2013).
\newblock Domain adaptation under target and conditional shift.
\newblock In:  {\em International Conference on Machine Learning\/}. PMLR. pp.\
   819--827.

\end{thebibliography}

	\begin{figure}[t!]
		\centering
		% First subfigure
		\begin{minipage}[t]{0.48\linewidth}
			\centering
			\includegraphics[width=\linewidth]{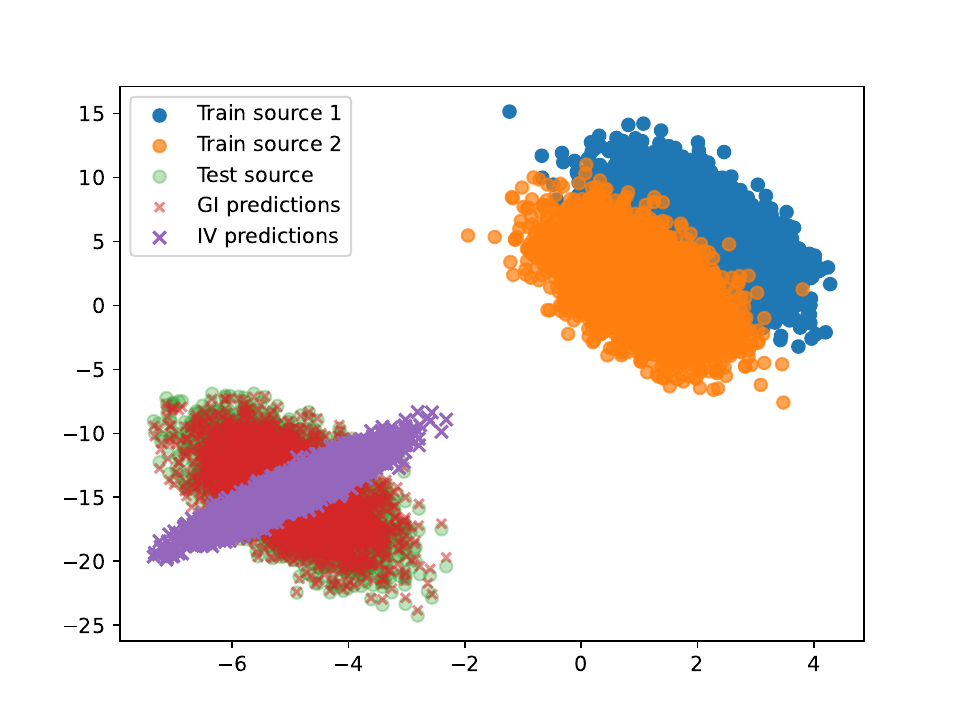}
			\caption{(Simpson's paradox) Although IV accurately identifies the causal parameter (positive slope), it proves entirely suboptimal for making predictions in the test environment. Our predictions, in red, align perfectly with the test source.}
			\label{fig:simpson}
		\end{minipage}%
		\hfill
		% Second subfigure
		\begin{minipage}[t]{0.48\linewidth}
			\centering
			\includegraphics[width=\linewidth]{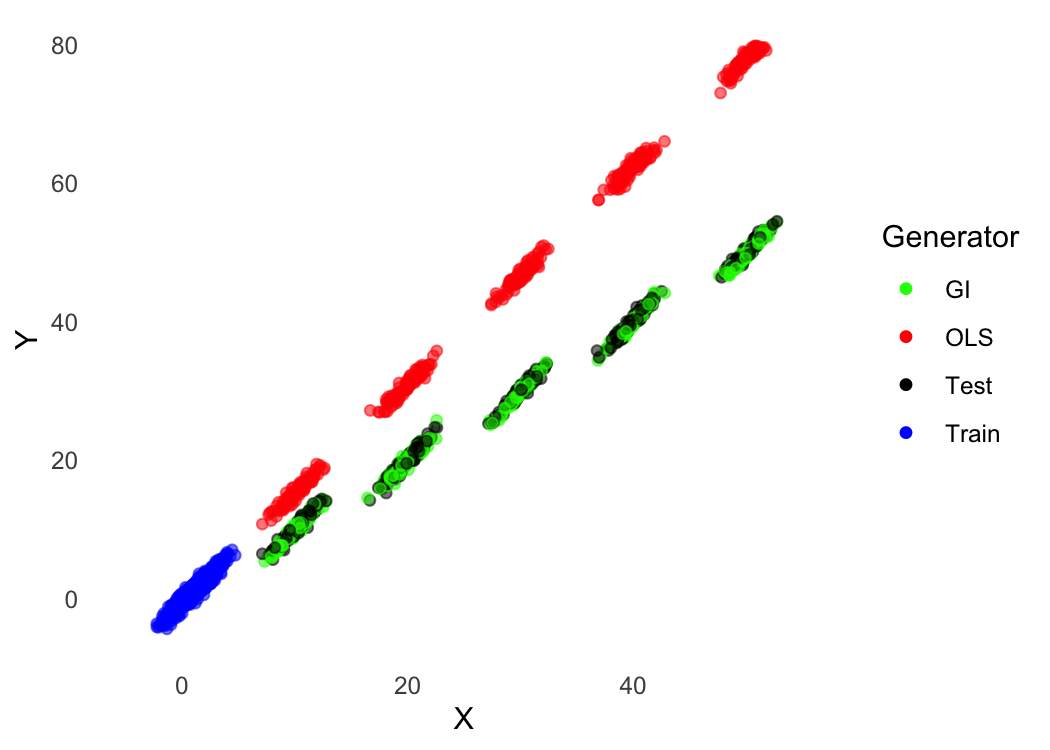}
			\caption{Sample from training environment $\mathbb{P}^1$ in blue. Samples from 5 different test environments in black. We fit an ordinary linear model to the blue data and launch red predictions. Our predictions, coloured in green, match the test distributions.}
			\label{firstfig}
		\end{minipage}
	\end{figure}

	\begin{figure}[h!]
		\centering
		\includegraphics[width=1\linewidth]{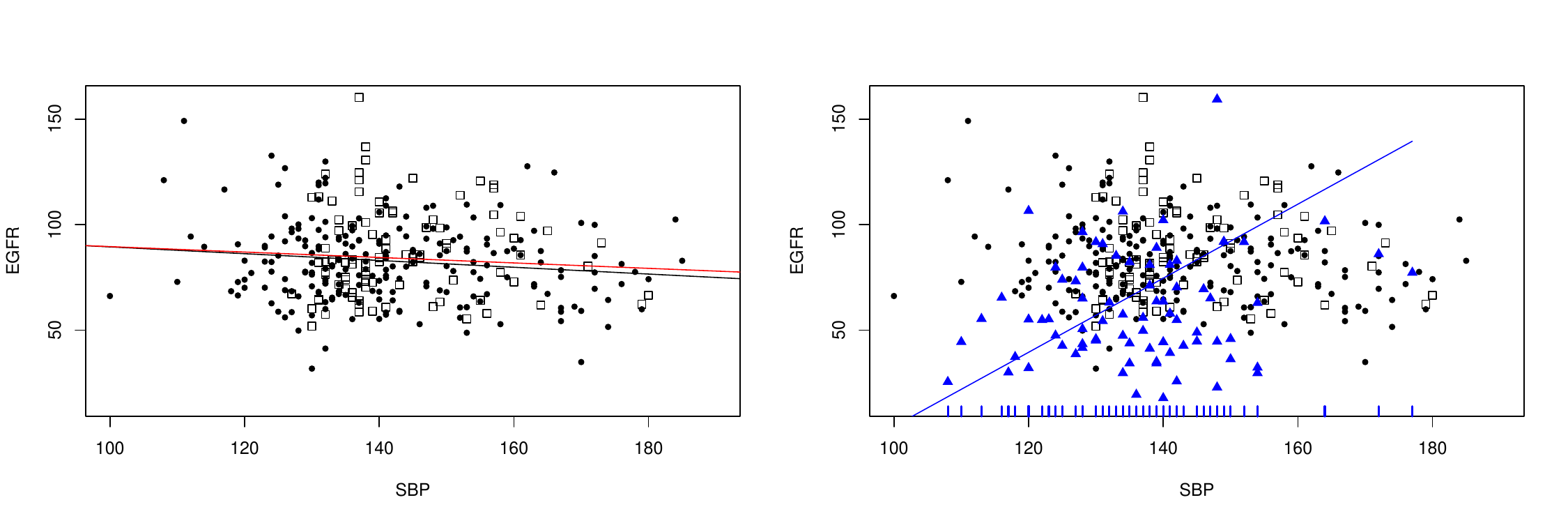}
		\caption{GI acting on SPRINT trial data. Left: data from hospital 48 in black dots and hospital 22 in white squares. Right: the predictions given by our estimator (in blue) for the new unseen hospital 15. Our approach does not see the blue triangles, just the blue values in the horizontal axis. }
		\label{in1}
	\end{figure}
	
		\begin{figure}[h!]
		\centering
		% First subfigure
		\begin{minipage}[t]{0.48\linewidth}
			\centering
			\includegraphics[width=\linewidth]{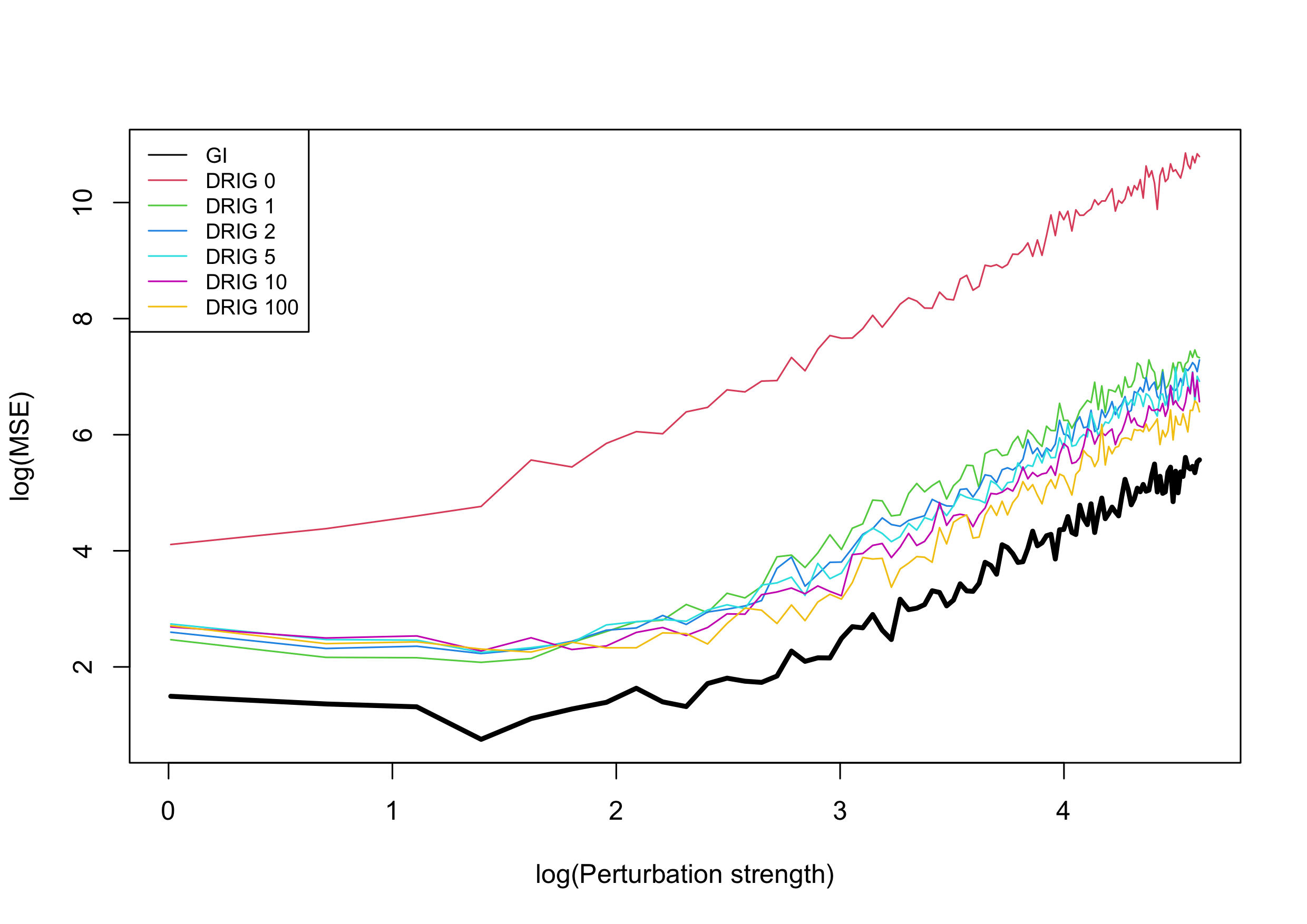}
			\caption{Test MSEs for varying test perturbation strengths. Notice the log scale in both axes, especially in the vertical one. GI's predictions MSE is always lower for any value of DRIG's tuning hyperparameter. }
			\label{simdrig}
		\end{minipage}\hfill
		% Second subfigure
		\begin{minipage}[t]{0.48\linewidth}
			\centering
			\includegraphics[width=\linewidth]{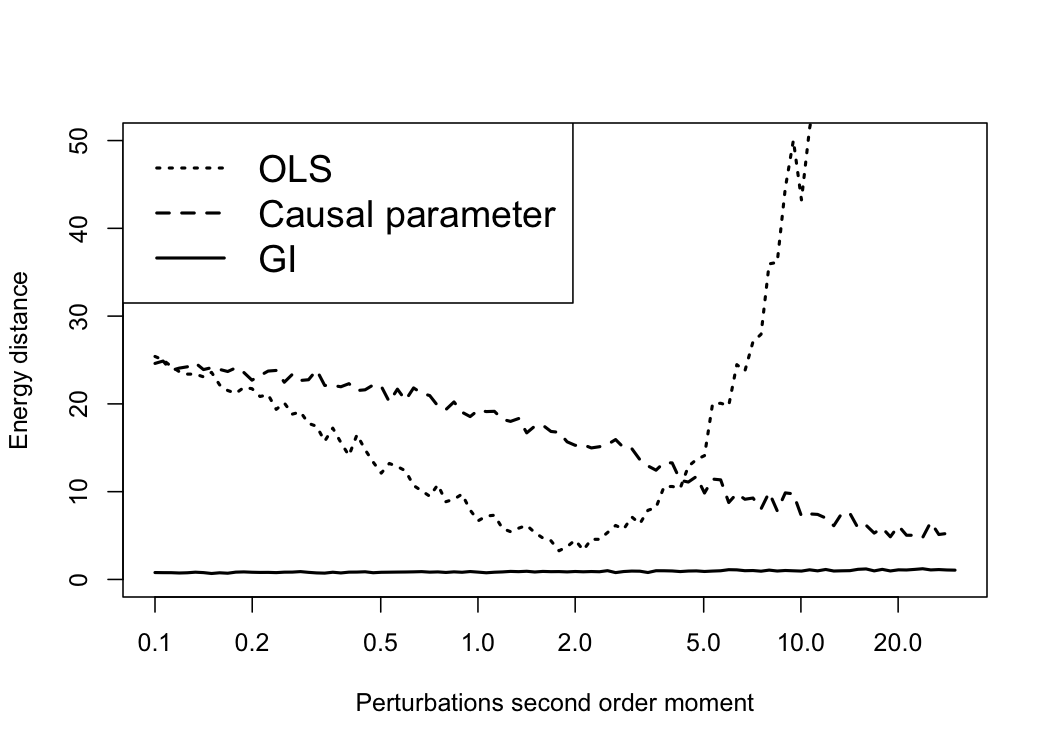}
			\caption{Performance of three generative models in emulating a test distribution indexed by the horizontal axis. Energy distance between GI's predictions and true test labels is uniformly almost zero.}
			\label{fig:energy}
		\end{minipage}

	\end{figure}
	
			\begin{figure}[t]
		\centering
		\includegraphics[width=1\linewidth]{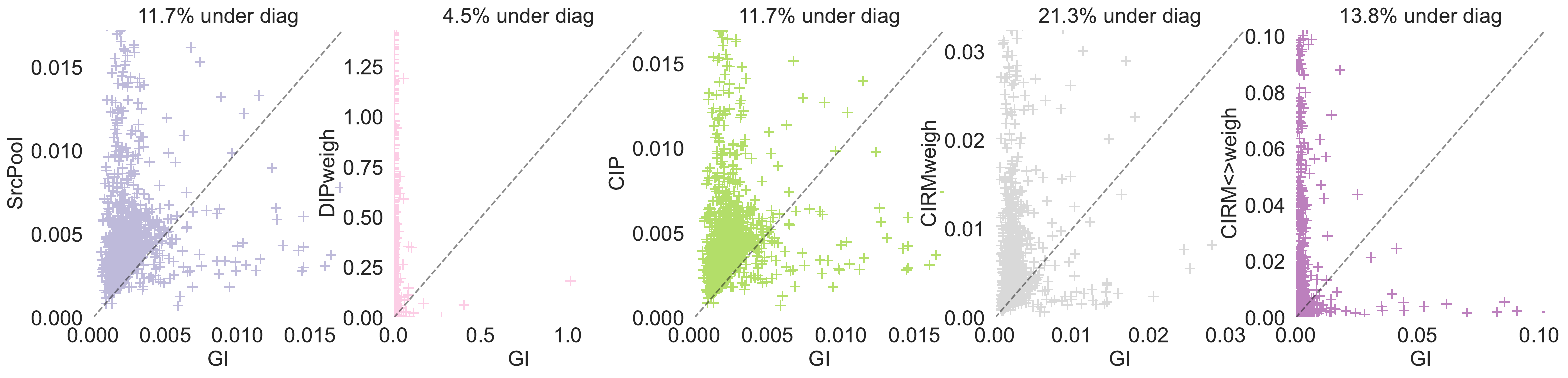}
		\caption{ Comparison of methods based on the number of times that GI outperforms in terms of MSE the other methods over 1000 random coefficient \textbf{B} data generations. }
		\label{fig:diag}
	\end{figure}
	
		\begin{table}[h!]
		\centering
		\caption{MSE summary for different sample sizes ($n$) with $1000$ simulation replications. It may be surprising that IV outperforms CD in a setting deemed favorable for CD according to its original paper. However, while this setting is advantageous for CD in terms of estimating $\beta_*$, it does not necessarily guarantee superior predictive performance.}
		\begin{tabular}{@{}cccccc@{}}
			\toprule
			\textbf{Sample Size} & \textbf{Method}       & \textbf{Median} & \textbf{25\% Quantile} & \textbf{75\% Quantile} \\ \midrule
			$n = 100$            & ${\text{IV}}$    & 14.1           & 7.0                   & 51.8                  \\
			
			& ${\text{CD}}$ & 34.5            & 10.5                  & 417.6                 \\
			& ${\text{GI}}$    & 12.6            & 6.6                   & 43.2                  \\  \midrule
			$n = 500$            & ${\text{IV}}$    & 7.3             & 5.5                   & 12.9                  \\
			
			& ${\text{CD}}$ & 12.8            & 7.0                   & 45.5                  \\ 
			& ${\text{GI}}$    & 6.8             & 5.5                   & 11.2                  \\ \midrule
			$n = 1000$           & ${\text{IV}}$    & 6.0             & 5.2                   & 8.6                   \\
			
			& ${\text{CD}}$ & 9.1             & 5.9                   & 22.4                  \\
			& ${\text{GI}}$    & 5.8             & 5.2                   & 7.9                   \\ \bottomrule
		\end{tabular}
		\label{simdantz}
	\end{table}

	\begin{table}[t!]
		\centering
		\begin{tabular}{|c|c|c|c|c|c|}
			\hline
			OLS & DRIG  $\gamma=0.05$ & DRIG $\gamma=0.1$ & DRIG $\gamma=0.3$ & DRIG $\gamma=0.5$ & \textbf{GI } \\
			\hline
			36.05 & 34.91 & 34.49 & 34.98 & 35.77 & \textbf{32.43} \\
			\hline
		\end{tabular}
		\caption{In each of the 10 test sources, MSE is computed, with the medians across the 10 sites displayed. }
		\label{numres}
	\end{table}
	
	\begin{figure}[t]
		\centering
		\includegraphics[width=0.9\linewidth]{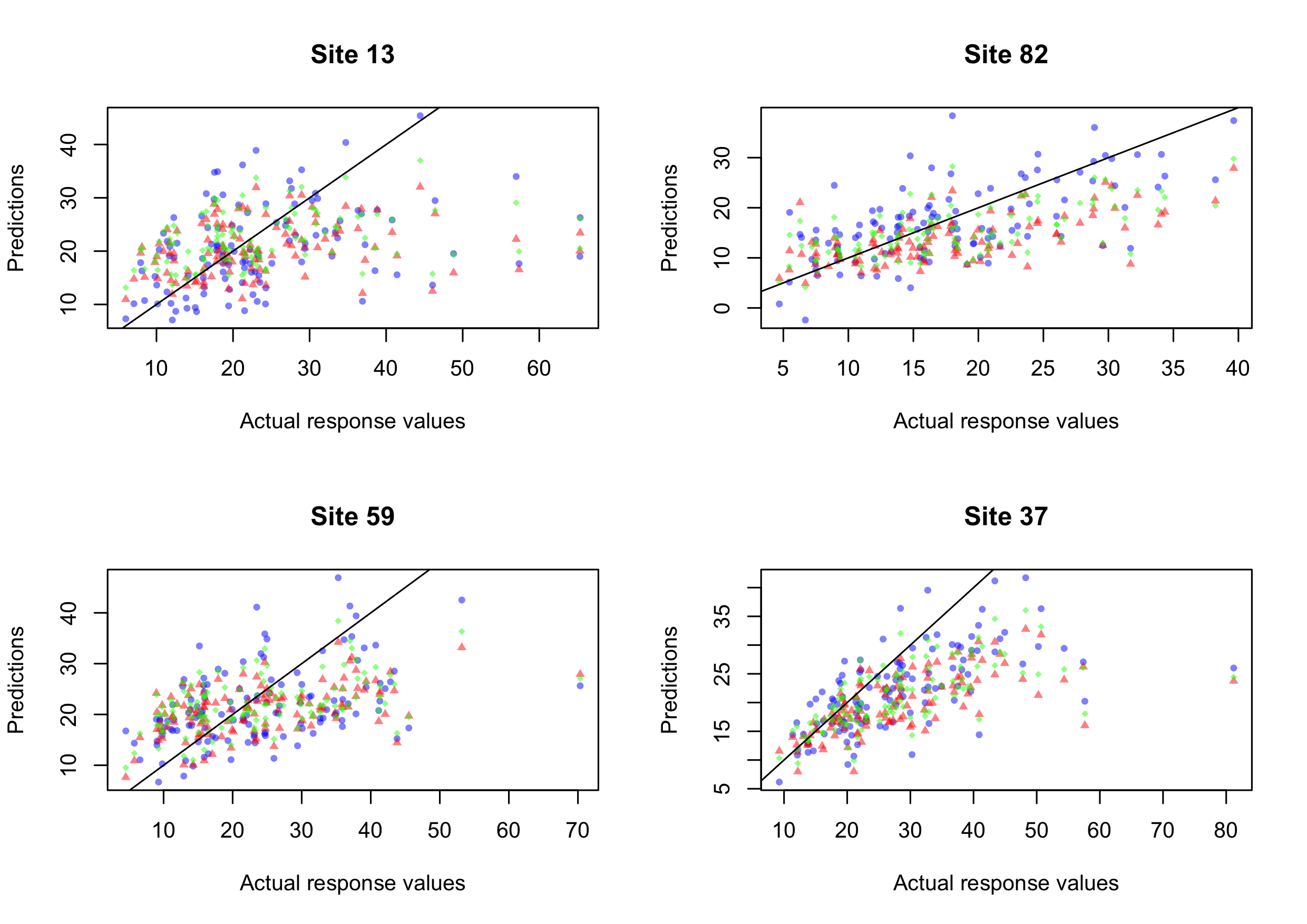}
		\caption{ Fitted vs. observed plots for four test sites. Green, red and blue correspond to OLS, DRIG $\gamma=0.1$ (best value) and GI, respectively. The solid black line is the identity function, where the points should ideally fall. }
		\label{testplots}
	\end{figure}
	\clearpage
	\appendix
	
	\section{Auxiliary results for asymptotic theory} 
	In this section we provide auxiliary results for asymptotic theory. We start by expressing the minimizer of theoretical risk (5) in terms of the observable random variables that represent each source $z=1,\ldots,Z$
	
	\begin{lemma}\label{popopt}
		$$\sum_{z=1} ^Z\mathbb{E}X_z\mathbb{E}X^T_z \beta_{opt} = \sum_{z=1} ^Z\mathbb{E}X_z\mathbb{E}Y_z $$ 
	\end{lemma}
	
	\begin{proof}	Recall that the theoretical risk is $\sum_{z=1}^Z\mathbb{E}\left[Y_z - \beta^TX_z - K^T(X_z - \mathbb{E}X_z ) \right] ^2$. Expanding the square inside of each expectation without substituting the model equation (1) and keeping the observable random variables
		$$
		\begin{aligned}  
			\sum_{z=1} ^Z\mathbb{E}[Y^2_z - 2\beta^TX_zY_z  - 2K^T(X_z - \mathbb{E}X_z )Y_z  +\\
			\beta^T X_zX^T_z\beta + 2\beta^T X_z(X_z-\mathbb{E}X_z)^TK +  \\
			K^T(X_z-\mathbb{E}X_z)(X_z-\mathbb{E}X_z)^TK].
		\end{aligned}
		$$
		Setting to zero the gradient wrt $\beta$ leads to
		$$
		\begin{aligned} \sum_{z=1}^Z\mathbb{E}X_zX^T_z \beta_{opt} + \mathbb{E}X_z(X_z-\mathbb{E}X_z)^TK = \sum_{z=1}^Z\mathbb{E}X_zY_z\\
		\end{aligned}
		$$
		Setting to zero the gradient wrt $K$ we get
		$$
		\begin{aligned} \sum_{z=1}^Z\mathbb{E}X_z(X_z-\mathbb{E}X_z)^T\beta_{opt} + \mathbb{E}X_z(X_z-\mathbb{E}X_z)^TK = \sum_{z=1}^Z\mathbb{E}Y_z(X_z-\mathbb{E}X_z)\\
		\end{aligned}
		$$
		Therefore
		$$\sum_{z=1}^Z\mathbb{E}X_z\mathbb{E}X^T_z \beta_{opt} = \sum_{z=1}^Z\mathbb{E}X_z\mathbb{E}Y_z.$$ 
	\end{proof}
	\noindent The next three results are the formal steps that lead to Corollary \ref{concmean}, which is a concentration inequality for the multivariate average of a random vector whose entries are sub-Gaussian. Lemma \ref{multimean} is a simple inequality holding with probability 1. 
	\begin{lemma}\label{multimean}
		Let ${X}_1, {X}_2, \ldots, {X}_n$ be i.i.d. samples from a $p$-dimensional random vector with mean ${\mu}$ and covariance matrix ${\Sigma}$ whose components are sub-Gaussian. Let the sample mean $\hat{\mu}$ be given by $\hat{\mu}=\frac{1}{n} \sum_{i=1}^n {X}_i$. Then we can write
		
		$$\|\hat{\mu}-{\mu}\|_2 \leq \sqrt{\frac{\lambda_{\max }({\Sigma})}{n}}\|W\|_2 \quad \textrm{a.s.}$$
		
		\noindent with $W$ a zero-mean $p$-dimensional random vector with identity covariance matrix whose components are sub-Gaussian. 
	\end{lemma}
	
	\begin{proof}
		We have that $\hat{\mu}-{\mu}$ has sub-Gaussian components as well (Lemma 2.6.8 in \cite{verhdp}) with mean zero and covariance matrix $\frac{1}{n} {\Sigma}$. Let the Cholesky decomposition of the covariance matrix be ${\Sigma}={L}{L}^T$;  where $L$ is a lower triangular matrix. We can write $
		\hat{\mu}-{\mu}=\frac{1}{\sqrt{n}} {L} {W}
		$
		where ${W} \sim \mathcal{N}\left(0, I_d\right)$. Taking the 2-norm of both sides:
		$\|\hat{\mu}-{\mu}\|_2=\left\|\frac{1}{\sqrt{n}} {L} {W}\right\|_2 
		=\frac{1}{\sqrt{n}}\|{L} {W}\|_2 $. Next,  $\|{L} {W}\|_2 \leq\|L\|_2\|{W}\|_2 $. The proof is completed by taking into account that $\|{L}\|_2=\sqrt{\lambda_{\max }({\Sigma})}$. 
	\end{proof}

	We use the following lemma to derive probabilistic bounds on the norm of sub-Gaussian random vectors. This lemma ensures that the norm of these vectors concentrates tightly around its expected value, which is crucial for establishing the convergence rates and concentration inequalities necessary for our asymptotic analysis. 
	
	\begin{lemma}(Concentration of the norm, Theorem 3.1.1 in \cite{verhdp}) Let $W=\left(W_1, \ldots, W_p\right) \in \mathbb{R}^p$ be a random vector with independent, sub-gaussian coordinates $W_i$ that satisfy $\mathbb{E} W_i^2=1$. Then

		$$\mathbb{P}\left\{\left|\frac{1}{p}\|W\|_2^2-1\right| \geq u\right\} \leq 2 \exp \left(-\frac{c p}{L^4} \min \left(u^2, u\right)\right) \quad \text { for all } u \geq 0$$
		
		$$
		\mathbb{P}\left\{\left|\|W\|_2-\sqrt{p}\right| \geq t\right\} \leq 2 \exp \left(-\frac{c t^2}{L^4}\right) \quad \text { for all } t \geq 0
		$$
		
		\noindent with $c$ an absolute constant derived by applying Bernstein's inequality and $L=\max _i\left\|W_i\right\|_{\psi_2}$, with $\|W\|_{\psi_2}=\inf \left\{t>0: \mathbb{E} \exp \left(W^2 / t^2\right) \leq 2\right\}$
	\end{lemma}
	
	\noindent Corollary \ref{corduo} is an alternative statement of the previous lemma. Its goal is to provide probability bounds that increase exponentially with the dimension \( p \) at the same rate for both the norm and its square, allowing us to combine both inequalities later on.
	\begin{Corollary}\label{corduo}
		Under the conditions of Lemma \ref{multimean} 
		
		$$\mathbb{P}\left\{|| W \|^2_2<p(L^4 + 1)\right\} \geq1 - 2 \exp \left(-cp \right)$$
		
		$$\mathbb{P}\left\{|| W \|_2<\sqrt{p}(L^2 + 1)\right\} \geq1 - 2 \exp \left(-cp \right)$$
	\end{Corollary}
	
	\begin{proof} For the first inequality, we impose $\frac{\min(u^2,u)}{L^4}=1$. Therefore, $\min(u^2,u)=L^4$ and this implies $u=L^4$ assuming $L\geq 1$, which makes sense according to proof of Theorem 3.1.1 in \cite{verhdp}. 
	\end{proof}
	
	\begin{Corollary}\label{concmean}
		Under the conditions of Lemma \ref{multimean}
		
		$$\mathbb{P}\left\{|| \widehat{\mu} - \mu \|_2>\sqrt{\frac{p\lambda_{\max }(\Sigma)}{n}}(L^2 + 1)\right\} \leq 2 \exp \left(-cp \right)$$
		
	\end{Corollary}
	
	\noindent Proposition \ref{univ} is a basic concentration inequality for the average of a one-dimensional sub-Gaussian random variable that can be found for example in \cite{verhdp}. 
	\begin  {Proposition}\label{univ} Let $Y_1,\ldots,Y_n$ be i.i.d. sub-Gaussian with norm $L$. Then,
		$$\mathbb{P}(|\overline{Y}-\mathbb{E}Y |\geq u) \leq 2 \exp \left(-\frac{n u^2}{2 L^2}\right)$$
	\end {Proposition}
	
	\noindent Define now
	
	$$\widehat{\Phi}= \frac{1}{Z}\sum_{z=1}^Z\widehat{\mu}_z\widehat{\mu}_z^T \quad \textrm{and} \quad \widehat{\Gamma}=  \frac{1}{Z}\sum_{z=1}^Z \widehat{\mu}_z \overline{Y}_z $$
	
	\noindent and 
	
	$${\Phi}= \frac{1}{Z}\sum_{z=1}^Z \mathbb{E}X_z\mathbb{E}X_z^T \quad \textrm{and} \quad \Gamma=\frac{1}{Z}\sum_{z=1}^Z\mathbb{E}X_z\mathbb{E}Y_z$$
	We need an intermediate concentration bound for the recently defined $\widehat{\Gamma}$. 
	
	\begin{Proposition} \label{gamma}	Let the assumptions in Lemma \ref{multimean} and Proposition \ref{univ} hold and let $$t:= Z\left[\sqrt{\frac{p}{n}}\left(\|\mu_M\|L\sqrt{2c} + \sqrt{\lambda_M}(L^2 + 1)|\mathbb{E}Y^M|\right) + \frac{p}{n}\sqrt{2\lambda_M}(L^2+1)L\right]$$

		$$\mathbb{P}\left(\|\widehat{\Gamma}-\Gamma\|_2 > t\right) \leq 8Z \exp{(-ct)} $$
	\end{Proposition}

	\begin{proof}
		$$ \widehat{\mu}\overline{Y} - \mu \mathbb{E}Y = \widehat{\mu}\overline{Y} - {\mu}\overline{Y} + {\mu}\overline{Y}  - \mu \mathbb{E}Y   $$
		
		\noindent Then, 
		
		$$
		\|\widehat{\mu}\overline{Y} - \mu \mathbb{E}Y\|_2 \leq \| \widehat{\mu} - {\mu} \|_2|\overline{Y}| + |\overline{Y}  - \mathbb{E}Y|\| \mu\|_2 \leq \| \widehat{\mu} - {\mu} \|_2|\overline{Y} - \mathbb{E}Y|+\| \widehat{\mu} - {\mu} \|_2|\mathbb{E}Y|  + |\overline{Y}  - \mathbb{E}Y|\| \mu\|_2 
		$$
		
		\noindent The proof is completed using the last inequality together with Corollary \ref{concmean} and Proposition \ref{univ}.
	\end{proof}
	
	Finally, we tackle each $z$ term in the definition of $\widehat{\Phi}$ and derive a probability bound for itself minus its respective population counterpart in $\Phi$.
	
	\begin{Proposition} \label{phi}
		Under the conditions of Lemma \ref{multimean},
		
		$$\left\|\hat{\mu} \hat{\mu}^T- {\mu}  {\mu}^T\right\|_2 = \mathcal{O}_{P}\left({\sqrt \frac{p}{n}}\right)$$
	\end{Proposition}
	
	\begin{proof}
		Start with the following decomposition
		$$\hat{\mu} \hat{\mu}^T- {\mu}  {\mu}^T=(\hat{\mu}- {\mu})(\hat{\mu}- {\mu})^T+ {\mu}(\hat{\mu}- {\mu})^T+(\hat{\mu}- {\mu})  {\mu}^T$$
		Therefore, 
		
		$$\left\|\hat{\mu} \hat{\mu}^T- {\mu}  {\mu}^T\right\|_2 \leq\left\|(\hat{\mu}- {\mu})(\hat{\mu}- {\mu})^T\right\|_2+2\left\| {\mu}\right\|_2\left\|\hat{\mu}- {\mu}\right\|_2$$
		On the one hand, 
		
		$$\left\|(\hat{\mu}- {\mu})(\hat{\mu}- {\mu})^T\right\|_2 = \lambda_{\max }\left((\hat{\mu}- {\mu})(\hat{\mu}- {\mu})^T\right)=\|\hat{\mu}- {\mu}\|_2^2$$
		\noindent	Finally, in virtue of both results in Corollary \ref{corduo}
		$$
		\left\|\hat{\mu} \hat{\mu}^T- {\mu}  {\mu}^T\right\|_2 \leq \frac{(L^4 + 1)\lambda_{\max }(\Sigma) p}{n} + 2(L^2 + 1)\left\| {\mu} \right\|_2 \sqrt{\frac{\lambda_{\max }(\Sigma) p}{n} },$$
		with probability greater than $1-4 \exp{(-cp)}$. 
		
		%$5$ stands for $L^4+1$ and $6$ stands for $2(L^2 +1)$
	\end{proof}

	\subsection{Proof of Theorem 5.1}

	From the normal equations and the population optimality conditions we have, respectively, $\widehat{\Phi} \widehat{\beta}= \widehat{\Gamma}$ and ${\Phi} {\beta}= {\Gamma}$. Then,

	$$\|\widehat{\Phi} - \Phi\|_2 \leq \frac{1}{Z}\sum_{z=1}^Z \|\widehat{\mu}_z \widehat{\mu}_z^T - \mathbb{E} X_z \mathbb{E} X_z^T \|,$$
	and
	
	$$\mathbb{P} \left( \|\widehat{\Phi}-\Phi\|_2 \leq \frac{1}{Z} \sum_{z=1}^Z \left(\frac{ \lambda_z p}{n}+2\left(L^2+1\right)\|\mu_z\|_2 \sqrt{\frac{\lambda_z p}{n}}\right)\right) \geq 1 - 4Z\exp{(-cp)}.$$

	\noindent If we take $$n > \frac{4}{\lambda^2_{\min}(\Phi)}\left((L^4+1) \lambda_M p + 2\left(L^2+1\right) \left\|\mu_M\right\|_2 \sqrt{\lambda_M p}\right)^2,$$ 
	then,
	
	$$\frac{\lambda_{\min}(\Phi)}{2} > \frac{1}{Z}\sum_{z=1}^Z \frac{(L^4+1) \lambda_z p}{n}+2\left(L^2+1\right)\left\|\mu_z\right\|_2 \sqrt{\frac{\lambda_z p}{n}}.$$
	Taking into account that $\Phi$ is symmetric and positive definite, we have by virtue of Corollary 6.3.4. in \cite{john} and taking, therein, $A=\Phi$ and $E=\widehat{\Phi} - \Phi$, that $\lambda_{\min}(\Phi) - \lambda_{\min}(\widehat{\Phi}) = |\lambda_{\min}(\Phi)| - |\lambda_{\min}(\widehat{\Phi}) | \leq |\lambda_{\min}(\Phi) - \lambda_{\min}(\widehat{\Phi}) | \leq \|\widehat{\Phi} - \Phi\|_2$. Therefore, 
	
	$$\mathbb{P} \left(\lambda_{\min}(\widehat{\Phi}) > \frac{\lambda_{\min}(\Phi)}{2}\right) \geq  1 - 4Z\exp{(-cp)}$$
	so invertibility of $\widehat{\Phi}$ happens with high probability. We can write 
	
	$$\widehat{\Phi} \widehat{\beta} -\widehat{\Phi} \beta_{opt}+ \widehat{\Phi} \beta_{opt}  - \Phi^* {\beta_{opt}} + \Phi^* {\beta_{opt}}=\widehat{\Gamma} - \Gamma + \Gamma.$$
	Equivalently, 
	
	$$
	\widehat{\Phi} (\widehat{\beta}- \beta_{opt}) + (\widehat{\Phi} - \Phi ^*)\beta_{opt} = \widehat{\Gamma} - \Gamma. 
	$$
	Therefore,
	$$
	\|\widehat{\beta}- \beta_{opt}\|_2 \leq \frac{2}{\lambda_{\min}(\Phi)} (\|\widehat{\Phi} - \Phi ^*\|_2 \|\beta_{opt}  \| _2+\|\widehat{\Gamma} - \Gamma \|_2 ).
	$$
	
	\noindent Finally, define $t$ in the statement of this Theorem as
	
	$$\begin{aligned}\frac{\lambda_{\min }(\Phi)}{2}t:=&\sqrt{\frac{p}{n}}\left(\left\|\mu_M\right\| L \sqrt{2 c}+\sqrt{\lambda_M}\left(L^2+1\right)\left|\mathbb{E} Y^M\right| +\left\|\beta_{opt }\right\|_2 2(L^2 + 1) \sqrt{\lambda_M} \right)+ \\&\frac{p}{n}\left( \sqrt{2 \lambda_M}\left(K^2+1\right) K + \left\|\beta_{opt }\right\|_2 (K^4 + 1) \lambda_M\right)\end{aligned}$$
	
	\noindent in order for it to be true. \hfill $\blacksquare$

	\subsection{Proof of Proposition 5.1}
	First, note that
	$$\begin{aligned}
		v^T\Omega_z v=&		(\mu_z^T\beta_*	)^2v^T\Sigma_z v- \mu_Z^T\beta_*v^T(K_* \mu^T_z +  \mu_z K_* ^T)v +v^T\mu_z\mu_z ^Tv\sigma^2 =\\
		=&	(\mu_z^T\beta_*	)^2v^T\Sigma_z v-2 \mu_Z^T\beta_*v^T K_*  \mu^T_zv +(v^T\mu_z)^2\sigma^2 \\
	\end{aligned}$$
	Using causal Mahalanobis condition $z$:
	$$v^T \Sigma_z v>\frac{1}{\sigma_Y^2}\left(v^T K\right)^2 \Longleftrightarrow\sigma_Y \sqrt{v^T \Sigma_z v} >\left|v^T K\right| >v^T K $$
	and, taking into account that there is a minus sign affecting the entire term, we substitute
	
	$$\begin{aligned}
		v^T	\Omega_z v>& 	(\mu_z^T\beta_*	)^2v^T\Sigma_z v-2 \mu_z^T\beta_*\sigma_Y \sqrt{v^T \Sigma_z v}   \mu^T_zv +(v^T\mu_z)^2\sigma^2 =\\
		=& 	\left(\mu_z^T\beta_* \sqrt{v^T \Sigma_z v}  - \mu^T_zv\sigma\right)^2 \geq 0 
	\end{aligned} $$ \hfill	$\blacksquare$

	\section{Proof of asymptotic normality}
	
	In this section we provide the proof of the asymptotic normality result (Theorem 5.2). In this proof, $z$ switches from a superscript to a subscript for the sake of writing convenience. Because of Propositions \ref{gamma} and \ref{phi} we have 
	
	$$(\hat{\Phi}, \hat{\Gamma})-(\Phi, \Gamma)=\mathcal{O}_P\left(\frac{1}{\sqrt{n}}\right).$$ 
	Note that 
	$$\begin{aligned}\widehat{\mu}_z\widehat{\mu}^T_z - \mathbb{E}\widehat{\mu}_z\widehat{\mu}^T_z = 
		\widehat{\mu}_z\widehat{\mu}^T_z - \mu_z \mu_z ^T +  \mu_z \mu_z ^T- \mathbb{E}\widehat{\mu}_z\widehat{\mu}^T_z = \mathcal{O}_P\left(\frac{1}{\sqrt{n}}\right) + \mathbb{E}\widehat{\mu}_z\mathbb{E}\widehat{\mu}^T_z - \mathbb{E}\widehat{\mu}_z\widehat{\mu}^T_z \\ = \mathcal{O}_P\left(\frac{1}{\sqrt{n}}\right) - \Cov \widehat{\mu}_z =  \mathcal{O}_P\left(\frac{1}{\sqrt{n}}\right) +  {o}_P\left(\frac{1}{{n}}\right). 
	\end{aligned}$$
	Proceeding analogously with $\widehat{\Gamma}$ we have
	
	$$(\hat{\Phi}, \hat{\Gamma})-(\mathbb{E}\hat{\Phi}, \mathbb{E}\hat{\Gamma})=\mathcal{O}_P\left(\frac{1}{\sqrt{n}}\right).$$ 
	Let us define the function $f: \mathrm{GL}_p \times \mathbb{R}^p \rightarrow \mathbb{R}^p$ by $f(X, Y):=X^{-1} Y$,
	where $X \in \mathrm{GL}_p$ (the set of invertible $p \times p$ matrices) and $Y \in \mathbb{R}$. From the proof of Theorem 1 in \cite{rothenhausler2019causal}, this function is continuously differentiable with derivative in direction $(D, d) \in \mathbb{R}^{p \times p} \times \mathbb{R}^p$ given by $ D f(X, Y)[(D, d)]=-X^{-1} D X^{-1} Y+X^{-1} d$. Therefore, by virtue of the delta method
	
	$$\begin{aligned}
		\widehat{\beta} - \beta _* = \hat{\Phi}^{-1}\hat{\Gamma} - {\Phi}^{-1}{\Gamma} = \left(-\Phi^{-1}(\hat{\Phi}-\Phi) \Phi^{-1} \Gamma+\Phi^{-1}(\hat{\Gamma}-\Gamma)\right)+{o}_P(1)  =  \\ \left(-\Phi^{-1}\left(\hat{\Phi}-\mathbb{E} \hat{\Phi} + {o}_P(1)\right) \Phi^{-1} \Gamma+\Phi^{-1}\left(\hat{\Gamma}- \mathbb{E} \hat{\Gamma} + {o}_P(1)\right)\right)+{o}_P(1).
	\end{aligned}
	$$
	Therefore,
	
	$$\begin{aligned}-\Phi^{-1}\hat \Phi \Phi^{-1}\Gamma+   \Phi^{-1} \widehat \Gamma= \frac{1}{Z}\sum_{z=1}^Z -\Phi^{-1} (\frac{1}{n} \sum_{i=1}^nX^z_iX_i^{z,T} - \widehat{\Sigma}_z) \Phi^{-1}\Gamma+   \Phi^{-1} \sum_{i=1}^n Y^z_i \widehat{\mu}_z=  \\
		\frac{1}{n}\sum_{i=1}^n\frac{1}{Z}\sum_{z=1}^Z\left\{  -\Phi^{-1} (X^z_i\mu^{z,T} + \mu{^z}X^{z,T}_i - \mu{^z}\mu^{z,T} + o_P(1) ) \Phi^{-1}\Gamma+   \Phi^{-1} Y^z_i (\mu_z + o_P(1))\right\}.
	\end{aligned}$$
	Let us focus on the $z$-th summand
	$$ -\Phi^{-1} \left(X^z_i\mu^{z,T} + \mu{^z}X^{z,T}_i - \mu{^z}\mu^{z,T}\right) \Phi^{-1} \Gamma+   \Phi^{-1} Y^z_i \mu_z.$$
	From population optimality, it turns out that
	$$Y^z = X^T \Phi ^{-1}\Gamma + \varepsilon_Y\Longrightarrow \mu^zY^z = \mu^zX^{z,T }\Phi ^{-1}\Gamma + \mu^z\varepsilon_Y. $$
	Replacing in the previous expression
	
	$$\begin{aligned} -\Phi^{-1} (X^z_i\mu^{z,T} + \mu{^z}X^{z,T}_i - \mu{^z}\mu^{z,T}) \Phi^{-1} \Gamma+   \Phi^{-1}(\mu^zX_i^{z,T }\Phi ^{-1}\Gamma + \mu^z\varepsilon_Y)  = \\
		-\Phi^{-1} (X^z_i\mu^{z,T}  - \mu{^z}\mu^{z,T}) \Phi^{-1} \Gamma+   \Phi^{-1} \mu^z\varepsilon_Y= \\
		-\Phi^{-1}( (X^z_i  - \mu{^z}) \mu^{z,T}\beta_*+  \mu^z\varepsilon_Y) =
		-\Phi^{-1}( (X^z_i  - \mu{^z}) \beta^{*,T} - \varepsilon_Y I_p)\mu^z. 
	\end{aligned}$$
	As the expectation of this last expression is zero, we compute its covariance matrix by taking the expectation of itself times itself transposed
	
	$$\begin{aligned}
		( (X^z_i  - \mu{^z}) \beta^{*,T}-  \varepsilon_YI_p)\mu_z\mu_z ^T( (X^z_i  - \mu{^z}) \beta^{*,T}-  \varepsilon_YI_p)^T = \\
		(\mu_Z^T\beta_*	)^2(X^z_i  - \mu{^z})	(X^z_i  - \mu{^z})^T - (X^z_i  - \mu{^z}) \beta^{*,T} \mu_ z \mu^T_z  \varepsilon_Y -   \varepsilon_Y\mu_z \mu^T_z\beta^{*}(X^z_i  - \mu{^z}) ^T +\varepsilon^2_Y\mu_z\mu_z ^T =\\
		(\mu_Z^T\beta_*	)^2(X^z_i  - \mu{^z})	(X^z_i  - \mu{^z})^T - \mu_Z^T\beta_*((X^z_i  - \mu{^z})  \mu^T_z  \varepsilon_Y +  \varepsilon_Y\mu_z (X^z_i  - \mu{^z}) ^T) +\varepsilon^2_Y\mu_z\mu_z ^T\\
	\end{aligned}$$
	Taking the expectation, we arrive to
	
	$$\begin{aligned}
		\Omega_z:=(\mu_z^T\beta_*	)^2\Sigma_z - \mu_Z^T\beta_*(K_* \mu^T_z +  \mu_z K_* ^T) +\mu_z\mu_z ^T\sigma^2. \\
	\end{aligned}$$
	If we assume the environments are independent, then variance becomes additive, so that the scaling factor in the CLT is
	
	$$	{\left(\frac{\frac{1}{Z^2}\sum_{z=1}^Z\Phi^{-1} \Omega_z \Phi^{-1}}{n}\right)^{-\frac{1}{2}}}.$$
	If each training source $z$ originates a different number of observations $n_z$,  the correct way of looking at the sums is as summing $Z$ independent averages with seldom additive contributions to the asymptotic covariance matrix: 
	
	$$\sum_{z=1}^Z  w_z\frac{1}{n_z}\sum_{i=1}^{n_z} (\cdots),$$ becoming equal to
	
	$${\left({\Phi^{-1} \sum_{z=1}^Z\frac{w_z^2}{n_z}\Omega_z \Phi^{-1}}\right)^{-\frac{1}{2}}}$$\qed

	Note that if all the sources had generated the same amount of data $n$, the asymptotic covariance matrix would be 
	$$\left(\frac{\frac{1}{Z^2} \sum_{z=1}^Z \Phi^{-1} \Omega_z \Phi^{-1}}{n}\right)$$ The factor $\frac{1}{Z^2}$ originates from the factor $\frac{1}{Z}$ present in $\frac{1}{nZ}M_Z^T M_Z$ defined as $\widehat{\Phi}$ during the proofs involved in the asymptotic theory of our estimator. One could think there is something wrong here, as an arbitrary power could be chosen since there is no law of large numbers across data sources: there is no average to pursue. However, everything behaves correctly in terms of scaling: if we had not incorporated the factor $\frac{1}{Z}$ in $\widehat{\Phi}$ then we would neither introduce it in ${\Phi}$. As a result, we are eventually using ${Z\Phi}$, whose inverse shows up twice and therefore the factor $\frac{1}{Z^2}$ would be present anyway.

	\section{Practical considerations from asymptotic normality} \label{pracons}
	
	We provide an initial approximation on how to use the asymptotic theory satisfied by Generative Invariance to choose the best data sources feeding our estimator in terms of well-conditioning of the asymptotic covariance matrix in case that we have access to $Z > p$ environments. In such situation we would have $\binom{Z}{p}$ different (dependent) estimators available. Suppose we fix $C =\binom{Z}{p}$ and estimate $\widehat{\beta}_1, \ldots, \widehat{\beta}_C$. Define \begin{equation}\label{bootav}
		\overline{\beta}= \frac{1}{C} \sum_{c=1} ^C \widehat{\beta}_c
	\end{equation}
	
	As a consequence of Proposition 5 in the main text, we have 
	
	\begin{Corollary} \label{propokey} Let $0 \neq v \in \mathbb{R}^p$ and  causal Mahalanobis condition $z$ hold for all $z=1, \ldots, Z$ . Then 
		$$	v^T\left(\sum_{z=1}^Z \frac{w_z^2}{n_z} \Omega_z \right )v > \sum_{z=1}^Z \frac{w_z^2}{n_z} \left(\mu_z^T \beta_* \sqrt{v^T \Sigma_z v}-\mu_z^T v \sigma\right)^2 $$
	\end{Corollary}
	
	We start by noting a necessary condition for a direction to belong to the null space of the asymptotic covariance matrix, as a consequence of Corollary \ref{propokey}.
	\begin{Corollary}\label{nec1}
		Let the assumptions in Theorem 6 hold. If $0 \neq v \in \mathbb{R}^p $ is in the null space of the asymptotic covariance matrix then $\left(\mu_z^T \beta_* \sqrt{v^T \Sigma_z v}-\mu_z^T v \sigma\right)=0$ for all $z=1,\ldots,Z$.
	\end{Corollary}
	
	\noindent Let us focus now on the quantity involved in the statement of the Corollary above. 
	
	\begin{definition} \label{defkey}Let  $ \overline{\beta}$ be estimated as in (\ref{bootav}). Define the \textit{efficiency key} as the function 
		
		$$ Q : r \in \mathbb{R}^p \longmapsto \sum_{z=1}^Z \frac{n_z}{(\sum_{z=1}^Z n_z)^2}\left(\hat{\mu}_z^T \overline{\beta} \sqrt{r^T \widehat{\Sigma}_z r-\widehat{\mu}_z^T r\widehat{\sigma}}\right)^2 \geq 0.$$
	\end{definition}

	\begin{algorithm}[ht!]\label{algo}
		\caption{Choosing the Best Data Sources for Asymptotic Efficiency. }\label{alg:best_sources}
		\begin{algorithmic}[1]
			\State \textbf{Input: } Data from $Z$ sources, $Z > p$
			
			\State Randomly split the set of sources into equal-sized $S_1$ and $S_2$. 
			\State Using $S_1$, compute estimators $\{\widehat{\beta}_c\}_{c \in \mathcal{C}(S_1, p)}$ 
			\State Estimate $\overline{\beta}$ from $S_1$:
			\begin{equation}
				\overline{\beta} = \frac{1}{|\mathcal{C}(S_1, p)|} \sum_{c \in \mathcal{C}(S_1, p)} \widehat{\beta}_c
			\end{equation}
			\For{each $b \in \mathcal{C}(S_2, p) $}
			\State Compute $\widehat{\beta}_b - \overline{\beta}$
			\State Complete $\{\widehat{\beta}_b - \overline{\beta}\}$ to an orthogonal basis of $\mathbb{R}^p$ using Gram-Schmidt
			\begin{equation}
				\langle \widehat{\beta}_b - \overline{\beta} \rangle \oplus \langle r^b_2, \ldots, r^b_p \rangle
			\end{equation}
			\State Evaluate the efficiency key $Q$ for $r^b_j$:
			\begin{equation}
				\widetilde{Q}_b = \min_{2 \leq j \leq p} Q(r^b_j)
			\end{equation}
			\EndFor
			\State Rank $\widetilde{Q}_b$ to quantify the efficiency of each combination of $p$ environments in $S_2$.
			\State \textbf{Output:} Optimal combination of $p$ environments in $S_2$ 
			
		\end{algorithmic}
	\end{algorithm}

	Let us derive a different necessary condition from Corollary \ref{nec1} for a vector to belong to the null space of a singular covariance matrix. A direction in the null space of the covariance matrix of a random variable represents a direction with zero variance. Consequently, sampling from the centered random variable results in realizations that fall in the orthogonal subspace to the null space, which is the image of the covariance matrix in \(\mathbb{R}^p\). Therefore, a vector in the null space of a covariance matrix must be orthogonal to the observations generated by a random variable with that covariance matrix. This implies that using an unseen set of environments $b=1,\ldots,B$ we can complete $\hat{\beta}_b - \overline{\beta}$ to an orthogonal basis $\{\hat{\beta}_b - \overline{\beta}, r^b_2, \ldots, r^b_p\}$ of $\mathbb{R}^p$ such that
	
	$$\langle \widehat{\beta}_b - \overline{\beta}\rangle \oplus \langle r^b_2, \ldots, r^b_p \rangle \quad b=1,\ldots, B.$$
	
	\noindent Finally, we define $$\widetilde{Q}_b = \operatorname{min}_{2 \leq j \leq p} Q(r^b_j) \quad b=1,\ldots, B$$
	
	\noindent whose ranking provides a procedure to quantify how good is a particular combination $c$ of $p$ environments in terms of invertibility of the asymptotic covariance matrix of the estimator $\beta_c$ produced by them. 
	
	Algorithm 1 presents a practical step-by-step user guide to the concepts discussed in this section. $\mathcal{C}(S, p)$ denotes there all the subsets with $p$ elements of a set $S$. In Subsection 5.3, we illustrate with a medical dataset how the discussed arguments lead to the optimal choice of hospitals for feeding Generative Invariance. As in practical applications it is often the case that $Z$ is not much bigger than $p$, the combinatorial number $\binom{Z}{p}$ is not very large and therefore the procedure described in Algorithm 1 is computationally fast.

	\section{Proofs regarding population Generative Invariance} \label{apppop}

	\subsection{Proof of Propostion 2}
	Consider the theoretical risk in (4), we have
	$$
	\begin{aligned}
		\mathbb{E} (Y_1 - \beta^TX_1 - K^T(X_1 - \mathbb{E}X_1))^2 = \mathbb{E}_{train} (Y - \beta^TX - K^T(X - \mathbb{E}_{train}X))^2. 
	\end{aligned}
	$$
	By expanding the square inside the expectation,
	$$
	\begin{aligned}
		&\mathbb{E}_{train} (Y - \beta^TX - K^T(X - \mathbb{E}_{train}X))^2 = \mathbb{E}_{train} (\beta_*^TX + \varepsilon_Y - \beta^TX - K^T(X - \mathbb{E}_{train}X))^2 \\  & \mathbb{E}_{train} ((\beta_* - \beta)^TX + \varepsilon_Y  - K^T(X - \mathbb{E}_{train}X))^2 = \\ & \mathbb{E}_{train} ((\beta_* - \beta)^T\varepsilon_X + \varepsilon_Y  - K^T(\varepsilon_X - \mathbb{E}_{train}X))^2 =
		\\ & \mathbb{E}_{train} ((\beta_* - \beta)^T\varepsilon_X)^2 + \mathbb{E}_{train} (\varepsilon_Y  - K^T(\varepsilon_X - \mathbb{E}_{train}X))^2 \\&+ 2\mathbb{E}_{train} ((\beta_* - \beta)^T\varepsilon_X) (\varepsilon_Y  - K^T(\varepsilon_X - \mathbb{E}_{train}X)  )= \\ & (\beta_* - \beta)^T (  \Sigma_1 + \mu \mu ^T )(\beta_* - \beta)+ \sigma^2_Y + K^T  \Sigma_1 K - 2K^TK_* + 2(\beta_* - \beta)^T(K_* -   \Sigma_1 K)
	\end{aligned}
	$$

	\noindent Set the gradient with respect to $\beta$ and $K$ equal to zero and we obtain
	
	$$
	\begin{aligned}\label{optbeta}
		(  \Sigma_1 + \mu \mu ^T )(\beta_* - \beta_{opt}) =   \Sigma_1K_{opt} -K_*
	\end{aligned}
	$$
	
	$$
	\begin{aligned}
		\Sigma_1K_{opt} - K_*  =   \Sigma_1(\beta_* - \beta_{opt}) 
	\end{aligned}
	$$\qed

	%
	%\subsection{Proof of Theorem \ref{t1}}
	%
	%Existence can be shown using the fact that the risk is coercitive under the Theorem assumptions together with its continuity. The latter can be proven by appling Lebesgue's dominated convergence theorem.
	%
	%From the previous proof we have
	%$$
	%\begin{aligned}
	%\Sigma^{-1}_X(  \Sigma_0K_{opt} - K_*) = (\beta_* - \beta_{opt}) 
	%\end{aligned}
	%$$
	%
	%And therefore,
	%
	%$$
	%\begin{aligned}
	%(  \Sigma_0 + \mu \mu ^T )\Sigma^{-1}_X(  \Sigma_0K_{opt} - K_*) =  \Sigma_0K_{opt}- K_*
	%\end{aligned}
	%$$
	%
	%As a consequence,
	%
	%$$
	%\begin{aligned}
	%((  \Sigma_0 + \mu \mu ^T )\Sigma^{-1}_X - I)(  \Sigma_0K_{opt} - K_*) = 0 \\
	%\end{aligned}
	%$$
	%
	%If $  \Sigma_0K_{opt} = K_*$ had not hold, then $((  \Sigma_0 + \mu \mu ^T )\Sigma^{-1}_X - I)=0$, but then we would had $(  \Sigma_0 + \mu \mu ^T )=  \Sigma_0$ and this would necessarily imply that $\mu$ is the zero of $\mathbb{R}^p$.
	%
	%
	%
	%As every covariance matrix is positive definite, we can use Sherman-Morrison's formula to invert $  \Sigma_0 + \mu \mu ^T$: 
	%
	%
	%$$\left(  \Sigma_0+\mu \mu^{\top}\right)^{-1}=  \Sigma_0^{-1}-\frac{  \Sigma_0^{-1} \mu \mu^{\top}   \Sigma_0^{-1}}{1+\mu^{\top}   \Sigma_0^{-1} \mu}$$
	%
	%And we deduce $\beta_{opt}=\beta_*$ from plugging $K_{opt}=K_*$
	
	\subsection{Proof of Proposition 3} \label{proofmulti}
	Consider the risk minimization problem (5)
	$$
	\begin{aligned}\argminD_{(\beta,K) \in \Theta} \sum_{z=1} ^Z\mathbb{E}\left[Y_z - \beta^TX_z - K^T(X_z - \mathbb{E}X_z ) \right] ^2.
	\end{aligned}
	$$
	We proceed by expanding the square inside each expectation
	
	$$
	\begin{aligned}\label{opt1}&\sum_{z=1} ^Z\mathbb{E}\left[Y_z - \beta^TX_z - K^T(X_z - \mathbb{E}X_z ) \right] ^2= \\ &\sum_{z=1} ^Z\mathbb{E}\left[ (\beta_* -\beta - K)^TX_z + \varepsilon + K^T\mu_z\right] ^2= \\
		&\sum_{z=1} ^Z\mathbb{E}\left[ (\beta_* -\beta - K)^TX_zX^T_z(\beta_* -\beta - K)+ 2(\beta_* -\beta - K)^TX_e\varepsilon\right] +\\ &\mathbb{E}\left[ 2(\beta_* -\beta - K)^TX_z\mu^T_z K +  K^T\mu_zX^T_z K + \varepsilon^2 \right]= \\
		&(\beta_* -\beta - K)^T\sum_{z=1} ^ZS_z(\beta_* -\beta - K)+ 2Z(\beta_* -\beta - K)^TK_*+\\ &2(\beta_* -\beta - K)^T\sum_{z=1} ^Z\mu_z\mu^T_zK +  K^T\sum_{z=1} ^Z\mu_z\mu^T_z K + Z \mathbb{E}\varepsilon^2 \\
	\end{aligned}
	$$
	\noindent Setting the gradient wrt $K$ to zero, leads to the linear system of equations
	$$
	\begin{aligned}\label{opt2}
		&-2\sum_{z=1} ^ZS_z(\beta_* -\beta - K) -2ZK_* - 4\sum_{z=1} ^Z\mu_z\mu^T_z K + 2\sum_{z=1} ^Z\mu_z\mu^T_z(\beta_* -\beta ) +2\sum_{z=1} ^Z\mu_z\mu^T_z K = 0 \\
		&-\sum_{z=1} ^Z\Sigma_z(\beta_* -\beta) + \sum_{z=1} ^Z\Sigma_zK - ZK_* = 0
	\end{aligned}
	$$
	
	\noindent Setting the gradient wrt $\beta$ to zero, we get the system
	$$
	\begin{aligned}
		&-2\sum_{z=1} ^ZS_z(\beta_* -\beta - K) -2ZK_* - 2\sum_{z=1} ^Z\mu_z\mu^T_z K =0 \\
		&-\sum_{z=1} ^ZS_z(\beta_* -\beta) + \sum_{z=1} ^Z\Sigma_zK -ZK_* =0
	\end{aligned}
	$$
	
	\noindent By combining both equations, if the minimizer exists it must satisfy,
	
	$$\sum_{z=1} ^ZS_z(\beta_* -\beta) = \sum_{z=1} ^Z\Sigma_z(\beta_* -\beta) \Longleftrightarrow   \sum_{z=1} ^Z\mu_z \mu^T_z(\beta_* -\beta) = 0 $$
	If the matrix
	$$ \sum_{z=1} ^Z\mu_z \mu^T_z$$
	is full rank, then if the minimizer of the theoretical risk $\beta_{opt}$ exists, it must be unique and equal to the causal parameter $\beta_*$. \qed

	\section{Auxiliary linear algebra results} \label{auxlinalg}
	
	We formulate two new linear algebra results together with their proofs.

	\begin{lemma}\label{aux1}
		Let $Z \leq p$. If $v_1, \ldots, v_Z \in \mathbb{R}^p$ are linearly independent, then $ \sum_{z=1}^Zv_z v^T_z$ has rank $Z$. 
		
	\end{lemma}
	
	\begin{proof}
		Define $A:=\sum_{z=1}^Zv_z v^T_z$. Extend the linearly independent set to a basis of $\mathbb{R}^p$ and form $A'$ as $A$, but with the completion vectors. 
		Assume wlog that $Z < p$. $A + A'$ has rank $p$ because if we had $(A+A')v = 0$ then $v$ would need to be orthogonal to all the elements in a basis and, therefore, $v=0$. Next, $p=rank(A+A') \leq rank(A) + rank(A') \leq rank(A) + p-Z$ so $rank(A) \geq Z$. If $Z=p$, proceed with $A'=0$.

	\end{proof}
	
	\noindent If a vector is in the column space of $\sum_{z=1}^Zv_z v^T_z$ there exists a $w \in \mathbb{R}^p$ such that we can write this vector as $\sum_{z=1}^Zv_z (v^T_zw)$, which is in the span of $v_1, \ldots, v_Z$. The question is if the reciprocal holds, i.e. whether every vector in the span belongs to the column space, since it could be that adding extra rank one matrices could decrease the rank. The next result shows that this is not the case. 
	
	\begin {Proposition}\label{corinv}
		If $v_1, \ldots , v_Z \in \mathbb{R}^p$ span $\mathbb{R}^p$, then $\sum_{z=1}^Zv_z v^T_z$ has full rank.  
	\end {Proposition}

	\begin{proof}If $Z=p+1$ and $v_1, \ldots , v_Z$ span $\mathbb{R}^p$, then there exist $p$ linearly independent vectors among the original $Z$. Assume wlog that they are the first $p$. Build $\sum_{z=1}^pv_z v^T_z$, which is invertible accoding to the previous Lemma 2.2. By Sherman-Morrison's formula \cite{john}, $\sum_{z=1}^pv_z v^T_z + v_Z v^T_Z$ is invertible iff $v_Z^{\top} A^{-1} v_Z \neq -1$ which never happens because $A$ is positive definite: $u^T\sum_{z=1}^pv_zv^T_z  u = \sum_{z=1}^p (v^T_z  u)^2 >0, u \neq 0$, which always holds. 
	\end{proof}
	
	\section{Estimator of response noise variance} \label{sigmay}
	
	In Section 6.5, we required an estimator for $\sigma_Y$ to empirically study the distributional behavior of our generative model (3). Looking at the proof of Propostion 3 for $Z=1$, if we assume in addition existence and uniqueness of minimisers $\beta=\beta_*$ and $K= K_{opt}$, we have
	
	$$\min _{(\beta, K) \in \Theta} \mathbb{E}\left[\left(Y^1-X^1 \beta-\left(X^1-\mathbb{E} X^1\right) K\right)^2\right] = \sigma^2_Y + K_{opt}^T  \Sigma_1 K_{opt} - 2K_{opt}^TK_* = \sigma^2_Y - K^{T}_*\Sigma^{-1}_XK_*$$
	Therefore, 
	
	$$\sigma^2_Y =\mathbb{E}\left[\left(Y^1-X^1 \beta_*-\left(X^1-\mathbb{E} X^1\right) \Sigma^{-1}_XK_*\right)^2\right] +  K^{T}_*\Sigma^{-1}_XK_*$$
	In the multi source case, 
	
	$$\begin{aligned}
		&\frac{1}{\sum_{z=1}^Zn_z}\sum_{z=1}^Zn_z\sigma_z^2\\&=\frac{1}{\sum_{z=1}^Zn_z}\sum_{z=1}^Zn_z\mathbb{E}\left[\left(Y_z-X_z \beta_*-\left(X_z-\mathbb{E} X_z\right) \Sigma_X^{-1} K_*\right)^2\right]+K^T_{opt} K_*
	\end{aligned}$$
	
	\noindent where the emprical counterpart to the first term on RHS is the minimum empirical risk divided by $	\frac{1}{\sum_{z=1}^Zn_z}$ \footnote{In \texttt{R}, just take the \texttt{mean()} of the component-wise squared $	\sum_{z=1}^Zn_z$-dimensional vector involved in the empirical risk}. Therefore, we are able to recover the pooled noise variance. 
	%Having access to an estimate $\widehat{\sigma^2_Y}$ of $\Var \varepsilon_Y$ is inmediate by using the empirical counterpart of the latter expression. This enables estimation of $S^p_K$ as a ``normalizing constant'' by noticing that $\Var g_K(X_0,\xi) =\Var (g_K(X_0,\xi) - \xi S^p_K + \xi S^p_K =\Var (g_K(X_0,\xi) - \xi S^p_K) +\Var( \xi S^p_K )  = \Var (g_K(X_0,\xi) - \xi S^p_K) +S^2_p = \Var \varepsilon_Y $ and therefore $S^p_K= \sqrt{\Var \varepsilon_Y - \Var (g_K(X_1,\xi) - \xi S_p) }$ where all the elements on the RHS are estimable. 

	\section{Literature review on unsupervised domain adaptation}
	
	{Domain Invariant Projection} (DIP) \citep{dip13} seeks to learn a shared intermediate subspace that bridges the source and target domains. It projects data into a low-dimensional latent space where the empirical distributions of source and target data are aligned by minimizing their Maximum Mean Discrepancy (MMD). This approach builds on {Transfer Component Analysis} (TCA) \citep{tca}, which measures distributional differences in a lower-dimensional space rather than directly in a Reproducing Kernel Hilbert Space (RKHS). Conditional Invariant Residual Matching (CIRM) combines DIP with conditional invariance penalty (CIP) \citep{chen}. Unlike DIP, which maps source and target covariates into a common subspace, CIP leverages label information from multiple source environments to identify conditionally invariant components. CIP relies exclusively on training data, while CIRM first employs CIP to compute a linear combination of these components as a proxy for $Y$. To address label distribution shifts, CIRM subsequently applies DIP-style matching to the residuals obtained from regressing $Y$ on its proxy.
	
	\cite{chen} have shown that under a linear SCM, DIP achieves low target error when the prediction problem is anticausal and label distribution is unperturbed. However, DIP underperforms compared to a source-only estimator in scenarios with label distribution perturbations or when the prediction problem is causal or mixed-causal-anticausal. To address these limitations, \cite{chen} introduced CIRM and its variants, which outperform DIP in mixed-causal-anticausal scenarios and those with label distribution shifts. Since CIRM essentially applies DIP to label-distribution-adjusted data, a weighted variant, \texttt{CIRMweigh}, was proposed, similar to the weighted version of DIP, \texttt{DIPweigh}, crafted for multi-source data.  A formal definition of \texttt{CIRMweigh} can be found in Appendix A.1 of \cite{chen}.

	\section{Simulation details for Section 6.2}
	
	We first set $p=5$, \(\beta_* = (1, 5, 12, 6, 7)\)- where the first entry corresponds to the intercept term- and \(K_* = (-2.7, 2.1, 1.0, -1.7)\).  The $4$ last entries of \(X_z\) follow a multivariate normal distribution \(\mathcal{N}_4(\mu_z, \Sigma_z)\), where \(\mu^j_z \sim \mathcal{N}(0, s_2^2)\) for \(j = 2, \ldots,5\); \(\Sigma_z = A_z^{T} A_z + I_4\); and $s_2$ as defined below. Each entry of \(A_z\) is sampled independently from \(\mathcal{N}(0, s_1^2)\), with $s_1$ specified below. The response noise variance is set to be \(K_*^{T}(\Sigma_z)^{-1}K_* + 1\). For each environment \(z = 1, \ldots, 6\), we generate \(n = 200\) samples (even though $5$ environments would suffice). We compute our estimator and DRIG for several values of its regularization parameter \(\gamma\). This completes the training phase. Next, we generate 100 artificial test environments with different intervention strengths. Let \(S\) take values on a grid from 1 to 100, and define \(s_1 \sim U(0, S)\) and \(s_2 = S - s_1\). We maintain the same \(\beta_*\) and \(K_*\) to respect invariance. By varying \(s_1\) and \(s_2\), we allow for wilder shifts without constraining the interventional distributions to specific parameters. This generates trajectories of test MSE across different strengths for both our approach and DRIG with \(\gamma \in \{0, 1, 2, 5, 10, 100\}\). When $\gamma=0$, the method corresponds to OLS on the so-called \textit{reference environment} as described in \citep{shen2023causality}. When $\gamma=1$, it corresponds to pooled OLS. For $\gamma=100$, the DRIG estimator should approximate $\beta_*$. It is important to recall that DRIG does not account for the particular impact of confounding in each source even though its population version identifies $\beta_*$ when $\gamma \rightarrow \infty$. Finally, we repeat the procedure across 100 different training draws and average the trajectories.

	\section{Simulation details for Section 6.5}
	First, we clarify how point estimates of the OLS parameter and the causal parameter \(\beta_*\) lead to generators used for our comparison. 
	
	\begin{itemize}
		\item The empirical OLS induces a generative model defined by $(x,\xi) \longmapsto x\widehat{\beta}_{\operatorname{OLS}} + \widehat{\sigma}^2 \xi,$
		where $\widehat{\sigma}^2 $ is the empirical squared residual standard error using $\widehat{\beta}_{\operatorname{OLS}}$.
		
		\item $\widehat{\beta}$ alone given by empirical Generative Invariance induces a generative model that mimics the do-interventional distribution \citep{Pearl_2009} of $Y$ under do$(X=x)$ $(x,\xi) \longmapsto x\widehat{\beta} + \widehat{\sigma_*}^2 \xi,$ where $\widehat{\sigma}^2_* $ is the empirical squared residual standard error using $\widehat{\beta}$.
		\item Generative Invariance $(\widehat{\beta},\widehat{K})$ induces the generative model in expression (3) of the main document, with an estimator of $\sigma^2_Y$ given by $\widehat{\sigma^2_Y}=\frac{1}{n}\|Y-X\widehat{\beta}-(X-M)\widehat{K}\|^2 + \frac{\widehat{K}^2}{ \widehat{\Var X}}.
		$
		
		A justified derivation can be found in Supplement 5.
		
	\end{itemize}
	
	We fix one sampled training dataset and estimate the generator parameters using it. Then, for each different test interventional distribution of the covariate, we compute energy distance between the true sample and three artificial ones corresponding to each generator. We consider a grid of possible second-order moments for $V$ ranging from $0$ to $30$. For each value $s$ in the grid, we sample the variance of $V$ from \texttt{runif(0,s)} and get the squared mean by subtracting this value from $s$. Then, we generate one test dataset using this distribution for the perturbations. 
	
	We hide the response values and produce predictions using our proposed generator besides the two defined above. Finally, we compute three energy distances (see Supplement 6): each one comparing the original test data with the sample synthesized by each generator. We do this for 50 times for each $s$ and average the energy scores.

	\section{Energy statistics} \label{energy}
	We need a measure of similarity between $(X_{0,1},Y_{0,1})\ldots,(X_{0,n}Y_{0,n})$- being the response values therein hidden to our estimator- and an artificial sample. Energy distance $\mathcal{E}(Q,R)$ of the distributions of two random variables $Q$ and $R$ taking values in the same metric space is a statistical measure used to determine the similarity between two distributions \cite{szekely2017energy}. Originating from the concept of Newton’s gravitational potential energy, energy statistics leverage distances between data points in metric spaces.
	
	As a consequence of Proposition 3 in \cite{szekely2017energy}, for $Q,R$ bivariate random variables and $\delta$ the Euclidean distance therein, we have that  $$
	\mathbb{E}(Q,R)= 2 \mathbb{E} \delta(Q, R)- \mathbb{E}  \delta\left(Q, Q^{\prime}\right)- \mathbb{E}  \delta\left(R, R^{\prime}\right)=0
	$$
	holds if and only if $Q$ and $R$ are identically distributed, where tilde stands for an iid copy. 
	For us, $Q=(X_0,Y_0)$ and $R=(X_0,\widetilde{Y}_0)$, where $\widetilde{Y}^0$ follows the distribution of $m_{\theta}(X_0,\xi)$ conditional to $X_0$ for some $\theta \in \Theta$. 
	
	Suppose that $Q_1, \ldots, Q_{n_1}$ and $R_1, \ldots, R_{n_2}$ are independent random samples from the distributions of $Q$ and $R$, respectively. The two-sample energy V-statistic corresponding to the energy distance $\mathcal{E}(Q, R)$ is
	$$
	\begin{aligned}
		\widehat{\mathcal{E}}_{n_1, n_2}= & \frac{2}{n_1 n_2} \sum_{i=1}^{n_1} \sum_{m=1}^{n_2}\left|Q_i-W_m\right| 
		& -\frac{1}{n_1^2} \sum_{i=1}^{n_1} \sum_{j=1}^{n_1}\left|Q_i-Z_j\right|-\frac{1}{n_2^2} \sum_{\ell=1}^{n_2} \sum_{m=1}^{n_2}\left|W_{\ell}-W_m\right| .
	\end{aligned}
	$$
	For us $n_1=n_2$ because both samples from $Z$ and $W$ involve the same test covariate observations.

\end{document}